%% file: cas-dc-template.tex
\documentclass[a4paper,fleqn]{cas-dc}

\usepackage[authoryear,longnamesfirst]{natbib}

\def\tsc#1{\csdef{#1}{\textsc{\lowercase{#1}}\xspace}}
\tsc{WGM}
\tsc{QE}

\usepackage{paralist} 
\usepackage{subcaption} 
\usepackage{enumitem} 

\newcommand{\kh}{_{h}^{k}}
\newcommand{\khp}{_{h}^{k+1}}

\newcommand{\khphalf}{_{h}^{k+1/2}}
\newcommand{\khmhalf}{_{h}^{k-1/2}}

\newcommand{\average}[1]{\langle{#1}\rangle}

\begin{document}
\let\WriteBookmarks\relax
\def\floatpagepagefraction{1}
\def\textpagefraction{.001}

\shorttitle{}    

\shortauthors{}  

\title[mode = title]{A Systematic Modeling Framework for Dynamic Simulation of Fixed-Bed Reactors}  

\tnotemark[0] 

\tnotetext[0]{This project has been funded by the MissionGreenFuels project DynFlex under The Innovation Fund Denmark no. 1150-00001B.} 

%

\author[0]{Marcus Johan Schytt}[
    orcid=0009-0001-8143-1878,
]


\ead{mschytt@dtu.dk}


\credit{Conceptualization, Methodology, Software, Writing - Original Draft, Visualization}

\affiliation[0]{organization={Department of Applied Mathematics and Computer Science},
            addressline={Technical University of Denmark}, 
            city={Kongens Lyngby},
            postcode={DK-2800}, 
            country={Denmark}}

\author[0]{John Bagterp Jørgensen}[
    orcid=0000-0001-9799-2808,
]
\cormark[0]

\ead{jbjo@dtu.dk}


\credit{Conceptualization, Methodology, Writing - Review \& Editing, Supervision, Funding acquisition}


\cortext[1]{Corresponding author}



\begin{abstract}
We present a modular and thermodynamically consistent modeling framework for simulating steady-state and transient behavior in fixed-bed reactors. Accurate simulation of dynamic reactor behavior is essential for enabling flexible operation in Power-to-X (P2X) applications, such as Power-to-Ammonia and Power-to-Methanol, where fluctuating renewable energy inputs demand robust and responsive process control. The proposed models integrate non-ideal thermodynamics through cubic equations of state and account for both advective and dispersive transport phenomena. We derive consistent mass and energy balances using internal energy as the energy state variable, and obtain temperature and pressure from thermodynamic constraints. Our simulation methodology provides the necessary model functions for steady-state and dynamic simulations, as well as parametric sensitivity analysis. It is applied to two fundamental fixed-bed reactor units, the fixed-bed reactor (FBR) and the direct-cooled reactor (DCR). In the context of ammonia synthesis, we simulate representative reactor variants, the adiabatic fixed-bed reactor (AFBR) and the isothermal direct-cooled reactor (IDCR). Simulations assess the impact of real and ideal thermodynamic models, transport assumptions, and steady-state approximations. Results show that real-fluid effects at elevated pressures significantly influence steady-state outlet temperatures and conversions for the IDCR, while common literature model assumptions generally provide accurate dynamic predictions. Altogether, the framework supports systematic reactor model development and analysis under variable operating conditions and model assumptions relevant to Power-to-X applications.
\end{abstract}

\begin{graphicalabstract}
\includegraphics[width=\linewidth]{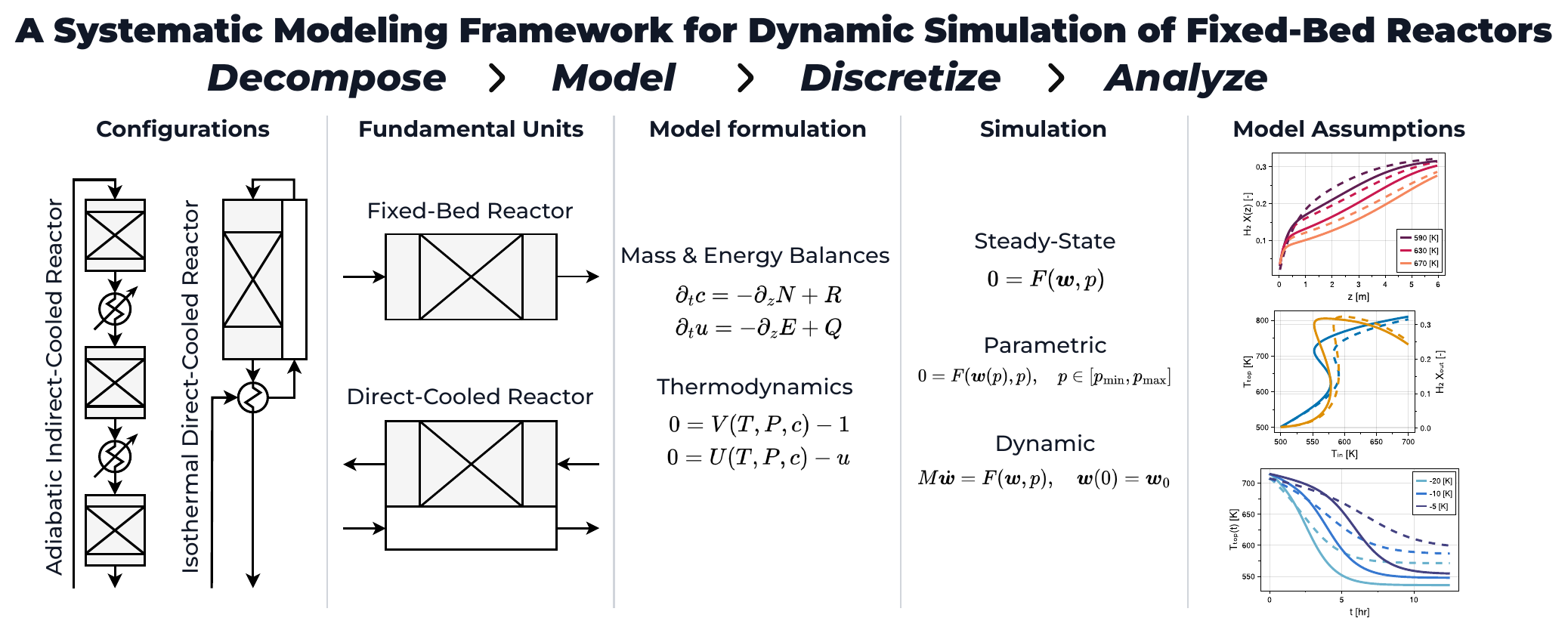}
\end{graphicalabstract}

\begin{highlights}
\item Modular modeling framework for fixed-bed reactor units with support for real-fluid thermodynamics.
\item Real-fluid thermodynamics are incorporated via cubic equations of state.
\item The framework enables steady-state and dynamic simulation using modern numerical methods.
\item Case study in ammonia synthesis shows the impact of reactor type and thermodynamic assumptions.
\item Simulation results verify common simplifying assumptions, but highlight the importance of real-fluid effects.
\end{highlights}


\begin{keywords}
Fixed-bed reactors\sep
Modeling framework\sep
Thermodynamics\sep
Partial differential-algebraic equations\sep
Steady-state simulation\sep
Dynamic simulation\sep
\end{keywords}

\maketitle

\input{tex/1introduction.tex}
\input{tex/2modeling.tex}
\input{tex/3simulation.tex}
\input{tex/4casestudy.tex}
\input{tex/5conclusion.tex}
\input{tex/6appendix.tex}










\printcredits

\bibliographystyle{cas-model2-names}

\bibliography{cas-refs}



\end{document}

%% file: tex/1introduction.tex
\section{Introduction}\label{sec:introduction}
\begin{figure*}[tb]
    \centering
    \includegraphics[width=\linewidth, trim={-1.5cm 0 0 0},clip]{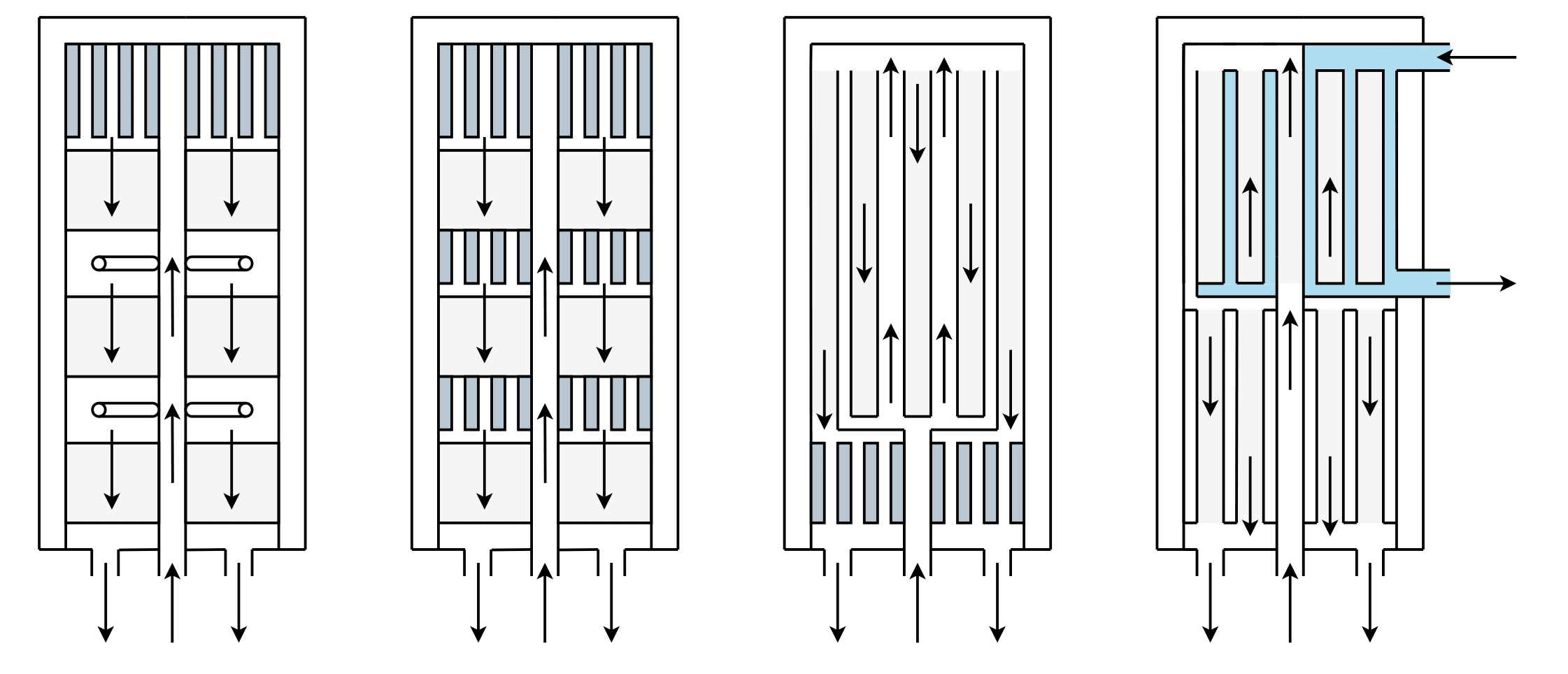}
    \caption{Schematic illustration of common commercial reactor configurations \citep{nielsen1995ammonia,bozzano2016efficient,khademi2017comparison}. \textbf{Left:} Adiabatic Quench-Cooled Reactor (AQCR). Comprises multiple adiabatic catalytic fixed beds with inter-bed cooling by injection of fresh feed gas. Industrial examples include the ICI low-pressure methanol converter and the M. W. Kellogg ammonia converter. \textbf{Center left:} Adiabatic Indirect-Cooled Reactor (AICR). Comprises multiple adiabatic catalytic fixed beds with inter-bed cooling via heat exchange with a cooling medium, often feed gas or water. Industrial examples include the Uhde ammonia converter and modern revamps of the M. W. Kellogg ammonia converter. \textbf{Center right:} Isothermal Direct-Cooled Reactor (IDCR). Comprises a multi-tubular shell-and-tube configuration of cooling tubes in a catalytic fixed bed, or of catalyst tubes immersed in a cooling medium, often feed gas or water. Industrial examples include the TVA ammonia converter (depicted here) and the Lurgi methanol converter. \textbf{Right:} Isothermal Hybrid-Cooled Reactor (IHCR). Comprises the connection of two IDCRs, one cooled by fresh feed gas and the other by water. Industrial examples include the Lurgi Megamethanol converter.}
    \label{fig:reactorconfigurations}
\end{figure*}

The European Union’s climate neutrality goals for 2050 \citep{EU2021delivering}, supported by national targets such as Denmark’s 70\% emissions reduction by 2030 \citep{Danish2020climate}, have enabled and accelerated the development of Power-to-X (P2X) technologies \citep{Danish2021strategy}. These technologies, which convert renewable electricity, typically wind and solar, into hydrogen via electrolysis. To enable large-scale storage, transport, and downstream utilization, hydrogen is further converted into energy-dense chemical carriers such as ammonia and methanol. These carriers are critical for decarbonizing hard-to-abate sectors, including heavy industry and long-distance transportation, such as aviation and shipping. However, the intermittent nature of renewable energy sources necessitates that P2X plants operate flexibly across a wide range of loads. Conventional reactor design and control strategies, which assume steady-state operation, are insufficient for such dynamic environments \citep{palys2022powertox}. Instead, growing emphasis is placed on design and control frameworks that optimize reactor performance under variable loads and enable real-time adaptation to fluctuating energy inputs \citep{rosbo2023flexible,rosbo2024optimal,fahr2024simultaneous,kong2024nonlinear, rosbo2026optimisation}. In this context, dynamic reactor models that accurately incorporate realistic thermodynamics and transport phenomena are indispensable tools for evaluating operational strategies, minimizing costs, and ensuring system stability in flexible P2X systems.

Modeling of fixed-bed reactors has been an active subfield of chemical engineering literature for over half a century. A major turning point came in the 1960s with the development of an analytical framework for transport processes, which laid the foundation for deriving reactor models grounded in the physical principles governing the conservation of momentum, mass, and energy. The seminal monograph by Bird, Stewart, and Lightfoot on transport phenomena~\citep{bird1960transport} has, over the decades, inspired numerous textbooks on the subject~\citep{froment1979chemical,rawlings2002chemical,jakobsen2014chemical}, ultimately shaping a mature and well-established field of reactor modeling. However, already by the early 1990s, it was recognized that effective equation-oriented numerical simulation required structured, high-fidelity process models \citep{marquardt1992objectoriented,pantelides1993equationoriented}.

Since it is neither practical nor necessary to model all physical phenomena in detail, the key scientific judgment lies in selecting assumptions appropriate to the specific reactor type and operating conditions. Conventional ammonia and methanol production operate under steady-state conditions. Accordingly, early reactor models~\citep{baddour1965steadystate,shah1967control,murase1970optimal} were formulated to simulate temperature profiles, identify hot spots, and perform parametric sensitivity analyses with respect to key operating parameters such as feed temperature, pressure drop, and catalyst activity. The analyses were used to investigate conditions leading to thermal runaway and to characterize instability phenomena such as extinction (blow-off). Beyond stability considerations, they also provided a basis for determining optimal operating conditions to maximize conversion and guide the design of reactor systems. At the time, however, extending these models to dynamic operation was considered computationally prohibitive~\citep{martinez1985modeling}.
Early dynamic models employed pseudo-steady-state assumptions, making only selected model equations dynamic. They were initially developed to study limit cycle instabilities~\citep{naess1992using,morud1998analysis, mancusi2009nonlinear} and the effects of catalyst degradation~\citep{lovik2001modelling}. Later contributions developed fully dynamic reactor models to investigate the transient and flexible operation of reactor systems~\citep{hansen1998plantwide,shahrokhi2005modeling,manenti2011dynamic,manenti2013dynamic, jorqueira2018modeling, rosbo2023flexible}.

To facilitate the modeling of fixed-bed reactors, various classifications have been proposed to systematically divide model types according to their spatial and phase resolution \citep{froment1979chemical}. Traditionally, the literature distinguishes between one- and two-dimensional models, depending on whether axial gradients alone are resolved or both axial and radial gradients are considered. In terms of phase resolution, models are typically categorized as pseudo-homogeneous, if no interfacial gradients between fluid and solid phases are modeled, or heterogeneous, if such gradients are explicitly accounted for. An alternative classification, proposed for the framework in~\citep{martinez1985modeling}, divides the reactor model into conceptual submodels, each treated separately with its own set of assumptions, allowing for modular model development. The constitutive equations typically involve a range of parameters, many of which are described using advanced empirical correlations from the literature~\citep{poling2000properties,bird2006transport,froment2010chemical}. The thermodynamic properties of are often simplified by assuming ideal fluid behavior. However, fixed-bed reactors typically operate at high pressures and over wide temperature ranges, where real-fluid effects can be significant, yet few studies account for them \citep{jorqueira2018modeling, rahmani2024design}.

Recently, \cite{martinsen2023modeling,rosbo2023flexible} introduced an approach to incorporate thermodynamic equations of state (EOS) into dynamic fixed-bed reactor models for the ammonia synthesis. In this work, we formalize and generalize this methodology by adopting a modular framework similar to~\cite{martinez1985modeling}. Specifically, we re-derive mass and energy balances using general transport assumptions from \cite{bird2006transport}, formulating the mass balances in terms of molar concentrations and energy balances in terms of internal energy density rather than temperature. The internal energy is explicitly described by the thermodynamic model, which simultaneously provides algebraic constraints from which reactor temperature and pressure profiles are computed. Unlike most literature models that assume ideal thermodynamics, our framework naturally accommodates real-fluid effects through the thermodynamic model. With this modeling framework, we aim to support systematic model development, steady and dynamic simulation, and flexible incorporation of model assumptions, adaptable across a variety of reactor configurations and operating conditions. We apply this framework to a set of fundamental reactor units represented in commercial reactor configurations for the ammonia and methanol synthesis (see Figure~\ref{fig:reactorconfigurations}). Specifically, we consider a case study in ammonia synthesis, where we describe variants of the reactor units and perform both steady-state and dynamic simulations to characterize steady-state curves, reactor profiles, and dynamic step responses. These simulations are used to explore optimal operating conditions under varying model assumptions.

\subsection{Overview}
We organize the rest of the paper as follows: Section~\ref{sec:modeling} introduces the modular modeling framework and presents the governing equations for fixed-bed reactor volumes. Section~\ref{sec:simulation} describes the numerical discretization using the method of lines and outlines the steady-state and dynamic simulation methodology. Section~\ref{sec:casestudy} applies the framework to a case study in ammonia synthesis and analyzes reactor units under different modeling assumptions. Section~\ref{sec:conclusion} summarizes the main findings and concludes the paper.

%% file: tex/2modeling.tex
\section{Modeling framework}\label{sec:modeling}

\begin{figure*}[htb]
    \centering
    \includegraphics[width=\linewidth]{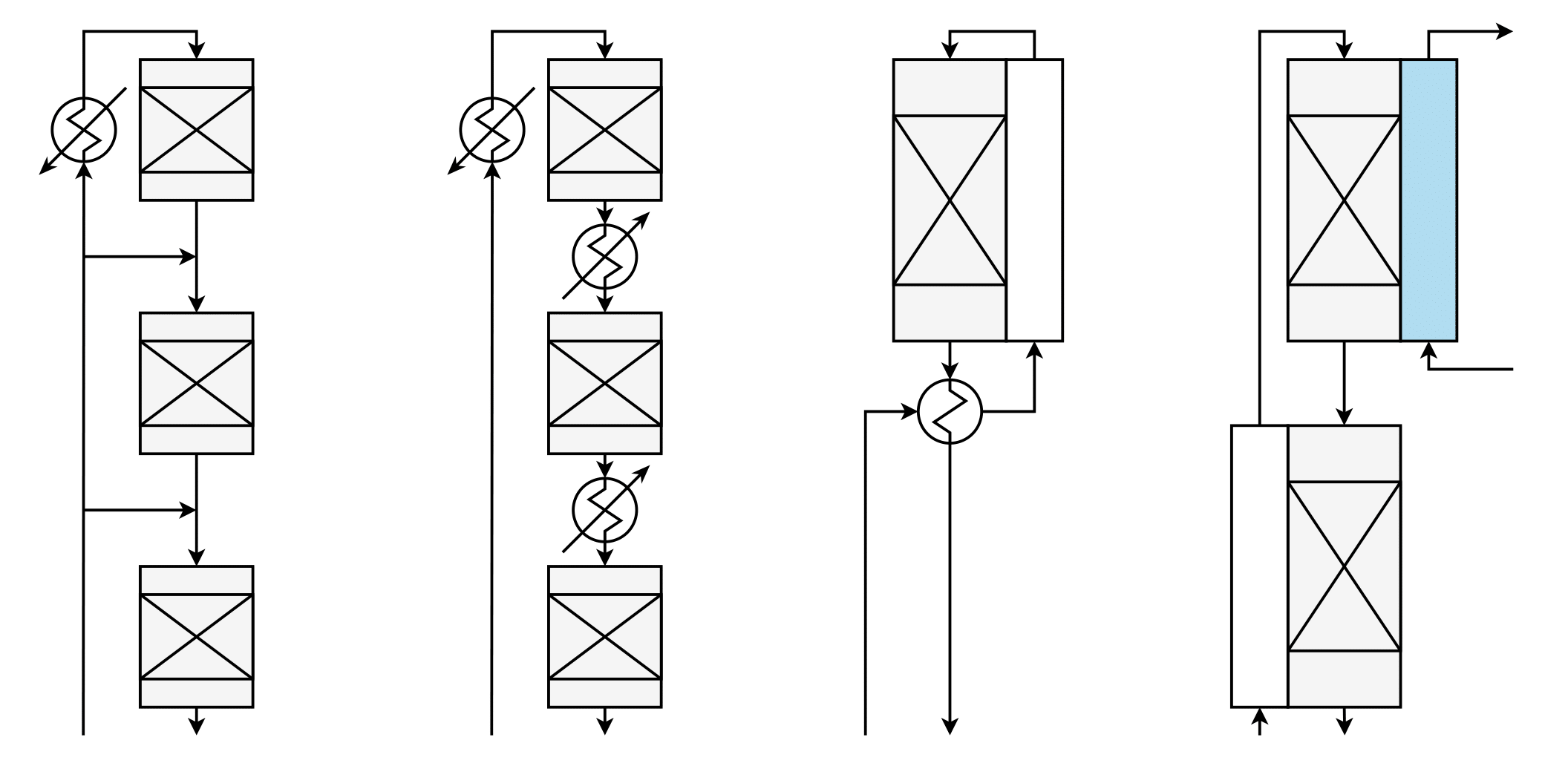}
    \caption{Diagrams of the reactor configurations in Figure~\ref{fig:reactorconfigurations}. \textbf{Left:} AQCR. \textbf{Center left:} AICR. \textbf{Center right:} IDCR. \textbf{Right:} IHCR.}
    \label{fig:reactordiagrams}
\end{figure*}

\begin{figure*}[htb]
    \centering
    \includegraphics[width=\linewidth]{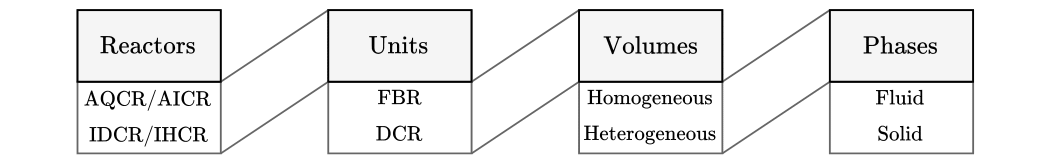}
    \caption{Diagram of the framework's model hierarchy. Each reactor configuration in Figure~\ref{fig:reactorconfigurations} is composed of reactor units, as shown in Figure~\ref{fig:reactordiagrams}, where each unit consists of a set of volumes that are either homogeneous fluid or heterogeneous fluid–solid.}
    \label{fig:modelhierarchy}
\end{figure*}

\begin{figure}[htb]
    \centering
    \includegraphics[trim={0 0.5cm 0 0.5cm},clip,width=\linewidth]{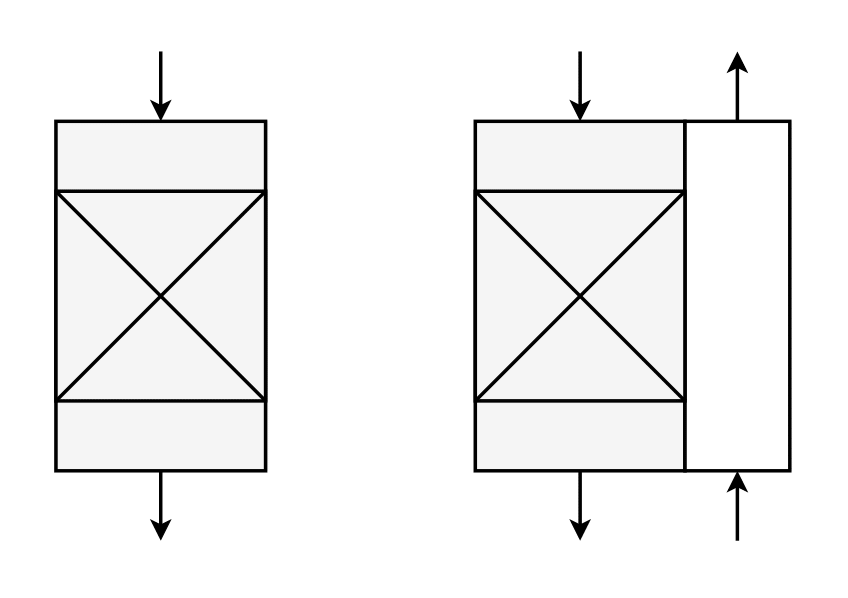}
    \caption{Diagrams of the units. \textbf{Left:} FBR. \textbf{Right:} DCR.}
    \label{fig:presentedunits}
\end{figure}

\begin{figure*}[htb]
    \centering
\includegraphics[width=\linewidth]{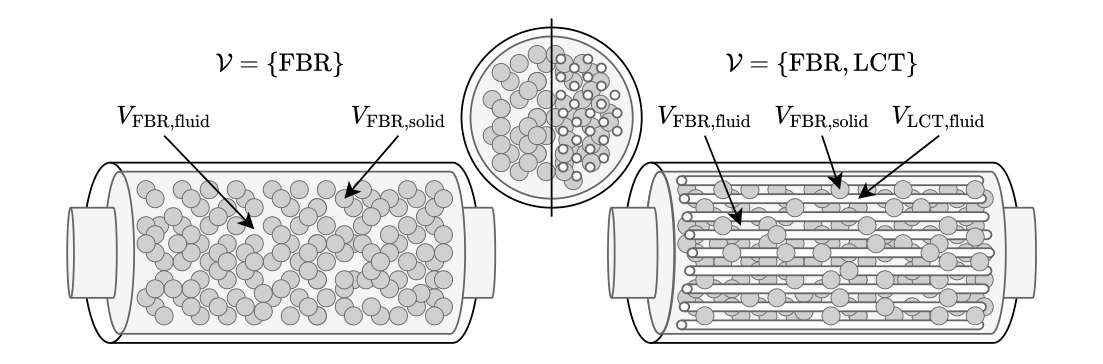}
    \caption{Schematic illustration of the units. \textbf{Left:} FBR. Modeled as one heterogeneous volume $V_\text{FBR}$ comprised of a fluid-phase volume~$V_\text{FBR,fluid}$ and a solid-phase volume~$V_\text{FBR,solid}$. \textbf{Right:} DCR. Modeled as two adjacent volumes, a heterogeneous fixed-bed reactor (FBR) volume $V_\text{FBR}$, and a homogeneous lumped cooling tube (LCT) volume~$V_\text{LCT}$. The LCT volume is entirely occupied by a fluid phase ($V_\text{LCT}=V_\text{LCT,fluid}$). \textbf{Center:} Cross-sectional view of the phase distribution inside the two respective reactor units.}
    \label{fig:reactorvolumes}
\end{figure*}

In this section, we present a general reactor modeling framework. We introduce a methodology for formulating systems of partial differential-algebraic equations (PDAEs) that describe fixed-bed reactors. As shown in Figure~\ref{fig:reactordiagrams}, each reactor is viewed as a collection of units that can be considered individually or as part of an interconnected system. Each reactor unit consists of volumes that are either homogeneous fluid or heterogeneous fluid–solid. Figure~\ref{fig:modelhierarchy} summarizes the reactor modeling hierarchy.

In this paper, we study the two reactor units in Figure~\ref{fig:presentedunits}:
\begin{itemize}
    \item Fixed-Bed Reactor (FBR).
    \item Direct-Cooled Reactor (DCR).
\end{itemize}
The FBR consists of one heterogeneous two-phase reactor volume. In contrast, the DCR is a multi-tubular shell-and-tube reactor volume with two subvolumes, a heterogeneous two-phase volume, and the other a homogeneous single-phase volume. The heterogeneous volume consists of a mobile fluid phase and a stationary solid catalyst phase, while the homogeneous volume consists of a single mobile fluid phase. For the DCR, the location of the heterogeneous subvolume, either on the shell-side or tube-side, is unimportant. In either case, the DCR can facilitate heat exchange via co- or counter-current flow in the homogeneous subvolume, possibly involving a different fluid. In this work, we will depict the DCR using a lumped cooling tube (LCT) with counter-current flow. Figure~\ref{fig:reactorvolumes} provides a graphical illustration of the two reactor units, including their subdivision into volumes and phases. We emphasize that the FBR and DCR units constitute a minimal yet sufficiently general set of building blocks, capable of describing all reactor configurations presented in Figure~\ref{fig:reactorconfigurations}.

We present first principles mass and energy balances, along with thermodynamic constraints, for each volume class. In describing these, the modeling framework adopts a modular structure, partitioning the volume model into five well-defined submodels:
\begin{enumerate}[label=\Roman*)]
    \item Thermodynamic model.
    \item Kinetic model.
    \item Mass transport model.
    \item Energy transport model.
    \item Heat transfer model.
\end{enumerate}
Each submodel is completed by a set of user-defined constitutive functions that encode thermodynamic properties, reaction kinetics, and transport relations. Separating physical laws from closure relations enables the independent specification or refinement of individual physical assumptions. Hence, the framework accommodates a wide range of reactor conditions and modeling fidelity without altering the overall framework structure. The modeling framework adopts a one-dimensional spatial resolution, wherein axial mass and energy propagation is denoted as transport, while interfacial exchange is referred to as transfer. We base our description of transport and transfer processes on the mechanistic principles presented in \cite{bird2006transport}. While empirical closures are not detailed here, relevant correlations can be obtained from established sources, such as the comprehensive compilation in \cite{poling2000properties}.

\subsection{Reactor volumes}
A reactor unit is partitioned into a set of cylindrical subvolumes $\mathcal{V}$. Each reactor volume is modeled as a cylindrical tube of length~$L$~[m] and volume~$V$~[m$^3$], assuming radial and azimuthal symmetry. Accordingly, the axial modeling domain is defined as the spatial interval $\Omega = [0, L]$.

\subsubsection{Geometric factors}
Transport and transfer processes within and between reactor volumes are described by several geometric factors.

\noindent\textbf{Heterogeneous volumes}\\
\noindent
 Heterogeneous volumes $V\in \mathcal{V}$ are partitioned into a fluid phase and a solid catalyst phase, such that $V = V_\text{fluid} + V_\text{solid}$. The fluid-phase volume fraction is defined as
 \begin{align}\label{eq:volumefraction}
     \varepsilon=\frac{V_\text{fluid}}{V}, \quad \frac{[\text{m}^3\text{-fluid}]}{[\text{m}^3]}\equiv[-],
 \end{align}
 hence the solid fraction is given by $1-\varepsilon$. For FBR volumes, $\varepsilon$ corresponds to the bed porosity, and we introduce the solid-to-fluid volume ratio
 \begin{align}\label{eq:volumeratio}
     \varphi=\frac{V_\text{solid}}{V_\text{fluid}}=\frac{1-\varepsilon}{\varepsilon}, \quad \frac{[\text{m}^3\text{-solid}]}{[\text{m}^3\text{-fluid}]}\equiv[-],
 \end{align}
 to describe catalytic reaction rates.
 
\noindent\textbf{Transport coupling}\\
\noindent When coupling axial fluid flow across adjacent volume boundaries, the differences in cross-sectional flow area must be accounted for. Under the assumed radial and azimuthal symmetry, the cross-sectional area of a volume $V \in \mathcal{V}$ is given by $S=V/L$ [m$^2$]. For a homogeneous volume, the fluid-phase area is simply $S_\text{fluid}=S$, while for a heterogeneous volume $S_{\text{fluid}}=\varepsilon S$. Given an adjacent volume $V' \in \mathcal{V}$, we define the fluid–to-fluid area ratio
\begin{align}\label{eq:fluidsratio}
    \psi = \frac{S'_{\text{fluid}}}{S_{\text{fluid}}} =\frac{V'_{\text{fluid}}}{V_{\text{fluid}}} , \quad [-],
\end{align}
to scale the interfacial fluxes.

\noindent\textbf{Interfacial transfer}\\
\noindent
The intensity of interfacial transfer processes depend on the extent of contact between volumes. If $A$ [m$^2$] denotes the interfacial surface area between two volumes $V,V' \in \mathcal{V}$, we quantify the interfacial contact per unit volume using the surface-area-to-volume ratio
\begin{align}\label{eq:surfacevolumeratio}
    a&=\frac{A}{V}, \quad [1/\text{m}],
\end{align}
and analogously for $V'$, to describe interfacial transfer rates.

\subsection{Balance equations and constraints}\label{sec:modelequations}
Each reactor volume $V\in\mathcal{V}$ is governed by its own particular set of mass and energy balances and thermodynamic constraints.

\subsubsection{Mass and energy balances}
\noindent\textbf{Mass balances}\\
\noindent
We study a fluid-phase mixture of chemical components~$\mathcal{C}$ and their molar concentrations $c=[c_\alpha]_{\alpha \in \mathcal{C}}$ [mol/m$^3$-fluid], defined per unit fluid volume. We assume that the components are subject to both advective and diffusive transport, as well as to reactive mass transfer. We summarize the mass balances by a system of transport equations
\begin{align}\label{eq:massbalances}
    \partial_tc=-\partial_zN+R.
\end{align}
The molar flux vector $N=[N_\alpha]_{\alpha \in \mathcal{C}}$ [mol$/$(s$\cdot$m$^2$-fluid)] is defined by the mass transport model, while the production term $R=[R_\alpha]_{\alpha \in \mathcal{C}}$~[mol$/$(s$\cdot$m$^3$-fluid)] is governed by the kinetic model, and expressed per unit fluid volume (in contrast to intrinsic rates defined per unit catalyst mass).

\noindent\textbf{Energy balance}\\
\noindent
The energy balance is written, neglecting kinetic and potential energy contributions, in terms of the internal energy density $u$ [J$/$m$^3$], defined per unit reactor volume. We assume that the internal energy is subject to convective and conductive transport, as well as interfacial heat transfer. In a similar fashion, we summarize the energy balance by a transport equation
\begin{align}\label{eq:energybalances}
    \partial_tu=-\partial_zE+Q.
\end{align}
Here we define the energy flux $E$~[W/m$^2$], which is defined by the energy transport model, and the heat transfer term $Q$~[W/m$^3$], which is defined by the heat transfer model.

\subsubsection{Thermodynamic constraints}
\noindent\textbf{Volume constraint}\\
\noindent The fluid volume is a thermodynamic property
\begin{align}
    V=V(T,P,n), \quad [\text{m}^3\text{-fluid}],
\end{align}
of the variables temperature $T$ [K], pressure $P$ [Pa], and component moles $n=[n_\alpha]_{\alpha \in \mathcal{C}}$ [mol]. Assuming local thermodynamic equilibria enables the pointwise use of classical thermodynamic relations \citep{jakobsen2014chemical}. As a consequence of the degree-one homogeneity of volume in its extensive variables \citep{michelsen2008thermodynamic}, and since concentrations are defined per unit fluid volume, we find the constant volume constraint
\begin{align}\label{eq:volumeconstraint}
    V(T,P,c)=1, \quad [-],
\end{align}
across the entire reactor length.

\noindent\textbf{Internal energy constraint}\\
\noindent The reactor volume's internal energy is likewise a thermodynamic property
\begin{align}
    U=U(T,P,n), \quad [\text{J}].
\end{align}
It gives rise to the internal energy density used in the energy balance \eqref{eq:energybalances}. Following similar arguments, we find the pointwise constraint
\begin{align}\label{eq:energyconstraint}
    U(T,P,c)=u, \quad [\text{J}/\text{m}^3],
\end{align}
which is elaborated below. Altogether, the mass and energy balances \eqref{eq:massbalances}-\eqref{eq:energybalances} provide closures for the differential variables $(c,u)$, whereas the thermodynamic constraints \eqref{eq:volumeconstraint}–\eqref{eq:energyconstraint} provide closures for the remaining algebraic variables $(T,P)$.

\subsection{Submodels}
\subsubsection{Thermodynamic model}\label{subsec:thermo}
We introduce a thermodynamic model to evaluate the thermodynamic properties required by the governing equations in Section~\ref{sec:modelequations}. Separate models are specified for fluid-phase and solid-phase properties. 

\noindent\textbf{Fluid model}\\
\noindent Given fluid-phase variables $(T,P,n)$ we require two functions to compute the thermodynamic fluid properties:
\begin{itemize}
    \item Volume $V=V(T,P,n), \quad$ [m$^3$-fluid].
    \item Enthalpy $H=H(T,P,n), \quad$ [J].
\end{itemize}
From these, the fluid-phase internal energy is recovered by the thermodynamic relation
\begin{align}\label{eq:thermorelation}
    U=H-PV,
\end{align}
We can compute the fluid-phase internal energy density directly from \eqref{eq:thermorelation} by defining the function
\begin{subequations}\label{eq:fluidenergydensity}
\begin{align}
    U(T,P,c)&=u_\text{fluid}(T,P,c)\\
    &=H(T,P,c)-PV(T,P,c).
\end{align}
\end{subequations}
The volume function is calculated using either the ideal gas law or extended to real fluids via a cubic equation of state (EOS) \citep{poling2000properties}. The enthalpy function is obtained by combining ideal heat capacity correlations with a residual correction derived from the EOS. Note that \eqref{eq:fluidenergydensity} simplifies when \eqref{eq:volumeconstraint} is satisfied. To facilitate numerical
simulation, as described in Section~\ref{sec:simulation}, we further
require the thermodynamic model's partial derivatives. They can be obtained either through analytical differentiation of the underlying thermodynamic models~\citep{ritschel2016opensource}, or by employing automatic differentiation frameworks~\citep{walker2022clapeyronjl,bell2022implementing}.

\noindent\textbf{Solid model}\\
\noindent In this work, we adopt the pseudo-homogeneous assumption, neglecting interfacial gradients and assuming that both phases share the same thermodynamic state variables $(T,P,c)$. Under this assumption, the volume’s internal energy density $u$ [J/m$^3$] is split into its phase contributions
\begin{align}
    u=\varepsilon u_\text{fluid}+(1-\varepsilon)u_\text{solid},
\end{align}
where $u_\text{fluid}$ [J/m$^3$-fluid] and $u_\text{solid}$ [J/m$^3$-solid] denote the internal energy densities of the fluid and solid phases, respectively. The scaling by the phase volume fractions $\varepsilon$ and $1-\varepsilon$ (as defined in \eqref{eq:volumefraction}), ensures that $u$ is defined per unit reactor volume. While $u_\text{fluid}$ is evaluated using \eqref{eq:fluidenergydensity}, we require a separate function
\begin{align}
    u_\text{solid}=u_\text{solid}(T,P),
\end{align}
to define one combined internal energy density function
\begin{align}
\begin{split}
    U(T,P,c)&=u(T,P,c)\\
    &=\varepsilon u_\text{fluid}(T,P,c)\\
    &\phantom{=}+(1-\varepsilon)u_\text{solid}(T,P),
\end{split}
\end{align}
for the entire reactor volume.

\subsubsection{Reaction model}\label{subsec:reaction}
\noindent\textbf{Stoichiometric model}\\
\noindent 
The production term $R$ [mol$/$(s$\cdot$m$^3$-fluid)] is written for a set of reactions $\mathcal{R}$ stoichiometrically so that the production rates are determined by a vector of reaction rates~$r=[r_k]_{k \in \mathcal{R}}$~{[mol$/$s$\cdot$m$^3$-fluid]} and by a constant stoichiometric matrix $\nu=[\nu_{k,\alpha}]_{(k,\alpha)\in\mathcal{R}\times\mathcal
C}$~[-]. The reaction rate of a component $\alpha \in \mathcal{C}$ is given by the stoichiometric relation
\begin{align}\label{eq:compreaction}
    R_\alpha = \sum_{k \in \mathcal{R}}\nu_{k,\alpha}r_k,
\end{align}
in which we do not discern between products and reactants. The combined production rate is then neatly summarized by
\begin{align}\label{eq:reaction}
    R=\nu^Tr.
\end{align}

\noindent\textbf{Kinetic model}\\
\noindent 
We assume that the reaction rates can be written as a function of the state variables
\begin{align}
    r=r(T,P,c),
\end{align}
so that the reaction rates may be described in terms of temperature, which is common for those including Arrhenius expressions, and pressure for those describing reactions between gasses. For non-ideal systems, the measure of component quantities may be better described by fugacities or chemical activities \citep{froment2010chemical}, which can be calculated by the thermodynamic model discussed above.

Heterogeneous catalytic reactions warrant a remark. For these, the above expressions \eqref{eq:compreaction}-\eqref{eq:reaction} remain valid, if the variables $(T,P,c)$ correspond to those at the catalytic zone. However, both interphase (external) and intraparticle (internal) mass transfer may severely restrict the apparent reaction rate, in which case the reaction is said to be diffusion-limited. Diffusion limitations have been studied for the ammonia synthesis in \cite{dyson1968kinetic, nielsen1995ammonia}, and for the methanol synthesis in \cite{lommerts2000mathematical}. External diffusion limitations can be accounted for by modeling interfacial mass transfer (adsorption), and similarly, internal diffusion limitations are most accurately modeled by including the internal particle mass balances in a heterogeneous model \citep{martinez1985modeling, lovik2001modelling, mancusi2009nonlinear}. However, explicitly modeling the intraparticle diffusive transport and interfacial transfer effectively makes the model two-dimensional and requires the specification of further state variables. Instead, we can compensate for diffusion limitations by including the effectiveness factors $\eta=[\eta_k]_{k\in\mathcal{R}}$ [-].
The effectiveness factors can be estimated from a generalized Thiele modulus \citep{lommerts2000mathematical, shahrokhi2005modeling}, or by literature correlations \citep{dyson1968kinetic}. For the sake of completeness, we assume them to be a function of the state variables
\begin{align}
    \eta = \eta(T,P,c).
\end{align}
Furthermore, to compensate for the fact that catalytic reaction rates are defined per unit solid volume, we scale the reaction rates by the volume ratio \eqref{eq:volumeratio}. Both scaling factors are taken into consideration, which leads us to the definition
\begin{align}\label{eq:pseudohomreaction}
r(T,P,c) \gets \frac{1-\varepsilon}{\varepsilon} \eta(T,P,c) \odot r(T,P,c).
\end{align}

\subsubsection{Mass transport model}\label{sec:masstransport}
The mass transport submodel is responsible for describing the molar flux vector $N$, which is split into its advective and diffusive parts
\begin{align}\label{eq:flux}
    N=N^\text{advection}+N^\text{diffusion}.
\end{align}

\noindent\textbf{Advection model}\\
\noindent
The advective flux $N^\text{advection}$ [mol$/$(s$\cdot$m$^2$-fluid)] describes the bulk flow of the fluid with component-wise velocities $v=[v_\alpha]_{\alpha \in \mathcal{C}}$ [m/s] following
\begin{align}\label{eq:advection}
    N^\text{advection}=v \odot c.
\end{align}
The velocities are usually set to a common scalar value which represents the molar average phase velocity. The velocity profile is subject to a momentum balance, which is pressure-driven and provides a semi-empirical relationship between the pressure gradient $\partial_z P$ [Pa/m] and the frictional forces, assuming negligible effects from gravitational forces, viscous forces, and convective acceleration~\citep{jakobsen2014chemical}. Correlations are usually parameterized by an appropriately defined Reynolds number
\begin{align}
    \text{Re}=\frac{\rho\ell|v|}{\mu},
\end{align}
in which $\rho$~[kg$/$m$^3$-fluid] denotes the fluid density, $\mu$ [Pa$\cdot$s] denotes the dynamic fluid viscosity, and $\ell$ [m] denotes some characteristic length scale. Following literature correlations, the fluid viscosity is described by a function of the state variables
\begin{align}
    \mu=\mu(T,P,c).
\end{align}
Moreover, the fluid density can be determined by the molecular weights $M=[M_\alpha]_{\alpha \in \mathcal{C}}$ [kg$/$mol], writing
\begin{align}\label{eq:fluiddensity}
    \rho=\rho(c)=\sum_{\alpha\in\mathcal{C}}M_\alpha c_\alpha.
\end{align}
If the chosen empirical law permits, entering it and the definition of the relevant Reynolds number into the corresponding pressure-drop relation will result in an explicit velocity equation
\begin{align}\label{eq:explicitvelocity}
    v=v(\partial_z P,\mu,\rho).
\end{align}

\noindent\textbf{Diffusion model}\\
\noindent
The diffusive flux $N_\text{diffusion}$~[kg$/$m$^3$-fluid] accounts for molecular and turbulent diffusion, and axial back-mixing caused by the presence of a stationary solid phase. Such effects are summarized by Fick's law of diffusion
\begin{align}\label{eq:diffusion}
    N^\text{diffusion}=-D\odot\partial_zc,
\end{align}
where $D=[D_\alpha]_{\alpha \in \mathcal{C}}$ [m$^2$$/$s] denotes the effective diffusion coefficients. The effective diffusion coefficients are typically assigned a single value for the mixture, which can be estimated using correlations from the literature. Again, we assume that the diffusion coefficients can be modeled as a function of the state variables
\begin{align}
    D=D(T,P,c).
\end{align}

\subsubsection{Energy transport model}\label{sec:heattransport}
Analogously, the Energy transport submodel is responsible for describing the the energy flux $E$, which is split into its convective and conductive parts
\begin{align}
    E=E^\text{convection}+E^\text{conduction}.
\end{align}

\noindent\textbf{Convection model}\\
\noindent
The convective flux $E^\text{convection}$ [W$/$m$^2$] describes the fluid-phase enthalpy flow. Consistent with prior assumptions, we neglect the contributions of kinetic and potential energy, as well as the work of viscous forces. The remaining fluid-phase contribution is neatly expressed using the advective molar flux and the degree-one homogeneity of the fluid-phase enthalpy property defined in Section \ref{subsec:thermo}. Accordingly, for homogeneous fluid volumes, we define
\begin{align}
 E^\text{convection}=H(T,P,N^\text{advection}), [\text{W/m}^2\text{-fluid}],
\end{align}
whereas for heterogeneous fluid–solid volumes under the pseudo-homogeneous assumption, we define
\begin{align}
    E^\text{convection}=\varepsilon H(T,P,N^\text{advection}), \quad [\text{W/m}^2].
\end{align}

\noindent\textbf{Conduction model}\\
\noindent
The conductive flux $E^\text{conduction}$ [W$/$m$^2$] accounts for diffusion-thermo effects and thermal conduction
\begin{align}
    E^\text{conduction}=q^\text{diffusion}+q^\text{thermal}.
\end{align}
The diffusion-thermo effect is defined by
\begin{align}
    q^\text{diffusion}=- \bar{H}^TN^\text{diffusion},
\end{align}
where $\bar{H}=[\bar{H}_\alpha]_{\alpha \in \mathcal{C}}$~[J$/$mol] denotes the partial molar enthalpies. Analogous to the molar diffusive flux, the thermal conductive flux is given by Fourier's law of thermal conduction
\begin{align}
    q^\text{thermal}=-\kappa\partial_zT,
\end{align}
where $\kappa$~[W$/$(m$\cdot$K)] denotes the effective conductivity. The partial molar enthalpies are defined by the EOS, and exploiting Euler's second theorem for homogeneous functions \citep{michelsen2008thermodynamic}, we write
\begin{align}
    \bar{H}=\partial_{c}H(T,P,c).
\end{align}
Moreover, the effective conductivity is determined from empirical correlations as a function of the states
\begin{align}
\kappa = \kappa(T, P, c).
\end{align}
For homogeneous fluid volumes, $\kappa$ corresponds to the effective fluid-phase conductivity, whereas for heterogeneous volumes under the pseudo-homogeneous assumption, $\kappa$ represents an effective combined conductivity
\begin{align}
\kappa = \varepsilon \kappa_\text{fluid} + (1 - \varepsilon) \kappa_\text{solid}.
\end{align}

\subsubsection{Heat transfer model}
The heat transfer term $Q$ [W$/$m$^3$] describes the interfacial transfer of energy across adjacent volumes due to heat conduction. Consider a neighboring volume $V' \in \mathcal{V}$ parallel to $V$, with state variables $(T', P', c')$. Newton's law of cooling defines the interfacial heat transfer by
\begin{align}
    Q=aK(T'-T),
\end{align}
where $a$ denotes the surface-area-to-volume ratio defined in~\eqref{eq:surfacevolumeratio}, and $K$ [W$/$(m$^2\cdot$K)] denotes the overall heat transfer coefficient. Note that $K$ captures the thermal resistance across the interface. For example, in the case of shell-and-tube configurations, $K$ can incorporate contributions from both fluid-side and shell-side heat transfer coefficients, as well as wall conduction resistance. Therefore it can be estimated by combining literature Nusselt number correlations using additivity of resistances. We define it as a function of the state variables of both volumes 
\begin{align}
    K=K(T',P',c'; T,P,c).
\end{align}

\subsection{Initial and boundary conditions}
The differential balances in \eqref{eq:massbalances} and \eqref{eq:energybalances} are completed by prescribing a set of initial conditions for the volume variables~$(c,u,T,P)$ across the entire reactor length $z~\in~\Omega$~[m] at time $t=0$~[s], and by specifying boundary conditions at the inlet $z=0$ and outlet $z=L$. We note that the present naming convention implies that flow occurs in the positive axial direction ($v>0$).

\subsubsection{Initial conditions}
We denote the initial conditions by
\begin{subequations}
\begin{align}\label{eq:ICs}
    c(0,z)&=c_0(z),\\
    u(0,z)&=u_0(z),\\
    T(0,z)&=T_0(z),\\
    P(0,z)&=P_0(z),
\end{align}
\end{subequations}
for all $z \in \Omega$. It is essential that the initial conditions are consistent with the thermodynamic constraints \eqref{eq:volumeconstraint} and \eqref{eq:energyconstraint} pointwise along the reactor.

\subsubsection{Boundary conditions}
We denote the boundary conditions (BCs) by
\begin{subequations}\label{eq:BCs}
\begin{align}
    N(t,0)&=N_\text{in}(t),\\
    E(t,0)&=E_\text{in}(t),\\ 
    N(t,L)&=N_\text{out}(t),\\
    E(t,L)&=E_\text{out}(t),
\end{align}  
\end{subequations}
for all $t \in [0,T]$, where $T$ [s] denotes the simulation horizon. We consider two types of BCs, as well as a set of coupling conditions.

\noindent\textbf{Flow-driven BCs}\\
\noindent
In accordance with the systematic flowsheeting procedure described in~\cite{hansen1998plantwide}, when the reactor volume is part of a reactor configuration or part of a larger flow sheet, the advective and convective flows are specified upstream, together with the downstream pressure. Hence, we specify:
\begin{itemize}
    \item Inlet molar flows $f_\text{in}=[f_{\text{in},\alpha}]_{\alpha \in \mathcal{C}}$ [mol/s].
    \item Inlet enthalpy flow $h_\text{in}$ [W].
    \item Outlet pressure $P_\text{out}$ [Pa].
\end{itemize}
Applying Danckwerts-type boundary conditions for the inflows, the corresponding inlet fluxes are then defined by
\begin{subequations}\label{eq:flowinlet}
\begin{align}
    N_\text{in}(t)&=\frac{f_\text{in}(t)}{S_\text{fluid}},\\
    E_\text{in}(t)&=\frac{h_\text{in}(t)}{S}.
\end{align}
\end{subequations}
Note that the inlet enthalpy flow need not be specified directly, as it can be computed from the upstream temperature and pressure $(T_\text{in},P_\text{in})$. Given the inlet molar flows $f_\text{in}$, the enthalpy flow can be evaluated by $h_\text{in} = H(T_\text{in}, P_\text{in}, f_\text{in})$,
using the fluid enthalpy function defined by the thermodynamic model. Similar to the inlet fluxes, the outlet fluxes are specified by the advective and convective flows. Hence, we obtain free outflow conditions by setting the diffusive and conductive fluxes to zero
\begin{subequations}
\begin{align}
    \partial_z c(t,L) &= 0,\\
    \partial_z T(t,L) &= 0.
\end{align}
\end{subequations}
The outlet molar fluxes are then given by the advective part
\begin{align}
    N_\text{out}(t)=v(\partial_zP(t,L),\mu(t,L),\rho(t,L))c(t,L).
\end{align}
Note that the outlet pressure $P_\text{out}(t)$ is specified to determine the outlet pressure gradient $\partial_z P(t,L)$. In a similar fashion, the outlet energy flux is given by the convective part. For homogeneous fluid volumes we define
\begin{align}
    E_\text{out}(t)=H(T(t,L),P(t,L),N_\text{out}(t)),
\end{align}
whereas for heterogeneous fluid-solid volumes, we define
\begin{align}
    E_\text{out}(t)=\varepsilon H(T(t,L),P(t,L),N_\text{out}(t)).
\end{align}

\noindent\textbf{Pressure-driven BCs}\\
\noindent
In contrast, when simulating a reactor volume in isolation, we may instead specify pressure-based boundary conditions. These do not rely on upstream flow data, but instead define the inlet state by prescribing the inlet pressure and composition. The following quantities are specified:
\begin{itemize}
\item Inlet molar fractions $x_\text{in}=[x_{\text{in},\alpha}]_{\alpha \in \mathcal{C}}$~[-].
\item Inlet temperature $T_\text{in}$~[K].
\item Inlet pressure $P_\text{in}$~[Pa]
\item Outlet pressure $P_\text{out}$~[Pa].
\end{itemize}
From these, the inlet fluxes are determined by
\begin{subequations}
\begin{align}
N_\text{in}(t) &= v_\text{in}(t) c_\text{in}(t), \\
E_\text{in}(t) &= H(T_\text{in}(t), P_\text{in}(t), N_\text{in}(t)),
\end{align}
\end{subequations}
where the inlet concentration $c_\text{in}=x_\text{in}/V(T_\text{in},P_\text{in},x_\text{in})$ is given by the thermodynamic fluid volume function, and the inlet velocity is reconstructed as
\begin{subequations}\label{eq:flowproperties}
\begin{align}
    v_\text{in}(t) &= v(\partial_z P(t,0),\mu_\text{in}(t),\rho_\text{in}(t)),\\
    \mu_\text{in}(t)&=\mu(T_\text{in}(t),P_\text{in}(t),c_\text{in}(t)),\\
    \rho_\text{in}(t)&=\rho(c_\text{in}(t)),
\end{align}
\end{subequations}
specified by the advection model. Here, the inlet pressure $P_\text{out}(t)$ is specified to determine the inlet pressure gradient $\partial_z P(t,0)$. The outflow conditions are the same as those for the flow-driven BCs.

\noindent\textbf{Coupling conditions}\\
\noindent
Finally, when a reactor unit consists of multiple adjacent volumes, internal coupling conditions must be imposed to ensure conservation of mass and energy across shared interfaces. Consider a neighboring upstream volume $V' \in \mathcal{V}$ with outlet fluxes $(N'_\text{out}(t),E'_\text{out}(t))$. To preserve continuity, the inlet fluxes of $V$ are defined by
\begin{subequations}\label{eq:couplingconditions}
\begin{align}
    N_{\text{in}}(t)&=\psi N'_{\text{out}}(t),\label{eq:convcoupling}\\
    E_{\text{in}}(t)&=\psi E'_{\text{out}}(t),
\end{align}
\end{subequations}
where $\psi$ is the fluid-to-fluid area ratio defined in \eqref{eq:fluidsratio}.

\subsection{Summary}
\begin{table}[]
\caption{Volume submodels and their constitutive functions.}
\centering
\resizebox{\columnwidth}{!}{%
\begin{tabular}{llr}
Model                            & Name                      & Symbol                    \\ \hline
\multicolumn{3}{l}{Thermodynamic model}                                                  \\ \hline
\multirow{2}{*}{Fluid model}     & Volume                    & $V(T,P,n)$                \\
                                 & Enthalpy                  & $H(T,P,n)$                \\
Solid model                      & Internal energy density   & $u_\text{solid}(T)$       \\ \hline
\multicolumn{3}{l}{Reaction model}                                                       \\ \hline
Stoichiometric model             & Stoichiometric matrix     & $\nu$                     \\
\multirow{2}{*}{Kinetic model}   & Reaction rates            & $r(T,P,c)$                \\
                                 & Effectiveness factor      & $\eta(T,P,c)$             \\ \hline
\multicolumn{3}{l}{Mass transport model}                                                 \\ \hline
\multirow{2}{*}{Advection model} & Viscosity                 & $\mu(T,P,c)$              \\
                                 & Velocity                  & $v(\partial_zP,\mu,\rho)$ \\
Diffusion model                  & Diffusivity               & $D(T,P,c)$                \\ \hline
\multicolumn{3}{l}{Energy transport model}                                                 \\ \hline
Conduction model                 & Conductivity              & $\kappa(T,P,c)$           \\ \hline
\multicolumn{3}{l}{Heat transfer model}                                                  \\ \hline
                                 & Heat transfer coefficient & $K(T',P',c';T,P,c)$       \\ \hline
\end{tabular}%
}
\label{tab:modelcomponents}
\end{table}

Each reactor volume is fully described by a set of submodel components, as outlined in Table~\ref{tab:modelcomponents}. Building on this general framework, we now specify the governing equations for the two general reactor units shown in Figure~\ref{fig:reactorvolumes}.

\subsubsection{Fixed-bed reactor}\label{subsec:FBR}
The FBR unit consists of a single heterogeneous reactor volume $\mathcal{V}=\{V_\text{FBR}\}$ with differential and algebraic state variables $(c,u,T,P)$. Its governing equations are
\begin{subequations}
\begin{align}
    \partial_tc&=-\partial_zN+R,\\
    \partial_tu&=-\partial_zE+Q,\\
    0&=V(T,P,c)-1,\\
    0&=u(T,P,c)-u,
\end{align}
\end{subequations}
with the production term defined by
\begin{align}
    R&=\nu^T r(T,P,c),
\end{align}
and the fluxes defined as
\begin{subequations}
\begin{align}
    N&=v(\partial_zP,\mu,\rho)\odot c-D(T,P,c)\odot\partial_zc,\\
    \begin{split}
    E&=\varepsilon \left(H(T,P,N^\text{advection})- \bar{H}^TN^\text{diffusion}\right)\\
    &\quad -\kappa(T,P,c)\partial_zT.
    \end{split}
\end{align}
\end{subequations}
In the absence of adjacent volumes, the heat transfer $Q$ term may be specified by an ambient temperature, or neglected altogether assuming adiabatic operation ($Q=0$).

\subsubsection{Direct-cooled reactor}\label{subsec:DCR}
The DCR unit is composed of two adjacent volumes $\mathcal{V}=\{V_\text{LCT},V_\text{FBR}\}$. We denote the state variables of the FBR by $(c,u,T,P)$, and those of the LCT by $(c',u',T',P')$. Their respective governing equations are summarized below.

\noindent\textbf{FBR subvolume}\\
The governing equations for the FBR subvolume match those of the standalone FBR unit, with the addition of the interfacial heat transfer term
\begin{align}
    Q=aK(T',P',c';T,P,c)(T'-T),
\end{align}
describing the coupling with the LCT.

\noindent\textbf{LCT subvolume}\\
Since no reaction occurs inside the homogeneous cooling tubes, its governing equations are summarized by
\begin{subequations}
\begin{align}
    \partial_tc'&=-\partial_zN',\\
    \partial_tu'&=-\partial_zE'+Q',\\
    0&=V(T',P',c')-1,\\
    0&=u_\text{fluid}(T,P,c)-u',
\end{align}
\end{subequations}
with fluxes defined by
\begin{subequations}
\begin{align}
    \begin{split}
    N'&=v'(\partial_zP',\mu',\rho')\odot c'\\
    &\phantom{=}-D'(T',P',c')\odot\partial_zc',
    \end{split}\\
    \begin{split}
    E'&=H(T',P',(N^\text{advection})')\\
    &\phantom{=}- (\bar{H}')^T(N^\text{diffusion})'\\
    &\phantom{=}-\kappa'(T',P',c')\partial_zT',
    \end{split}
\end{align}
\end{subequations}
and heat transfer term defined by
\begin{align}
    Q'=a'K(T',P',c';T,P,c)(T-T').
\end{align}
describing its heat exchange with the FBR volume.


\subsubsection{Literature model assumptions}\label{subsec:literaturemodels}
The governing equations presented here differ fundamentally from conventional reactor models commonly found in the literature. In Appendix \ref{app:energybalance}, we derive the conventional mass and energy balance equations by making explicit the key assumptions typically employed in the literature. In contrast, our formulation is based on fewer assumptions and more naturally accommodates real-fluid thermodynamics.

%% file: tex/3simulation.tex
\section{Simulation methodology}\label{sec:simulation}
To gain insight into the reactor's behavior under various operating conditions, we perform numerical simulations of the governing equations. The presented simulation methodology serves three primary objectives:
\begin{itemize}
    \item Determine steady-state profiles.
    \item Perform parametric sensitivity analysis.
    \item Simulate dynamic trajectories.
\end{itemize}
To this end, we apply the method of lines and first present the numerical discretization of the system of PDAEs governing the reactor volumes described in Section \ref{sec:modeling}.

\begin{figure}[tb]
    \centering
    \includegraphics[width=\linewidth]{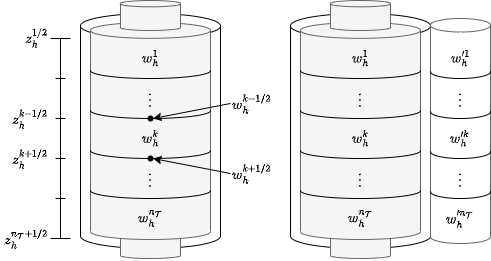}
    \caption{Schematic illustration of the finite volume scheme. \textbf{Left:} AFBR. \textbf{Right:} DCR. Divided into FBR and DCR volumes.}
    \label{fig:discretizedunits}
\end{figure}

\subsection{Spatial discretization}\label{subsec:spatial}
The spatial discretization is performed using a first-order upwind finite volume scheme \citep{leveque2002finite,hundsdorfer2003numerical}. We discretize the axial domain $\Omega$ uniformly into $n_\mathcal{T}$ finite volumes $\mathcal{T}\kh = [z\khmhalf, z\khphalf]$ defined by a characteristic length $h>0$, such that $\Omega = \cup_{k=1}^{n_\mathcal{T}}\mathcal{T}\kh$.
For each finite volume $\mathcal{T}\kh$, we define its midpoint and interfaces
\begin{subequations}
\begin{align}
    z\kh &= (k-1/2)h,\\
    z\khmhalf &= (k-1)h,\\
    z\khphalf &= kh,
\end{align}    
\end{subequations}
for $k=1,\dots,n_\mathcal{T}$. This is presented for both units in Figure~\ref{fig:discretizedunits}. Given a state variable $w(t,z)$ on $[0,T]\times\Omega$, we study its finite volume discretization $w_h=[w\kh]_{k=1}^{n_\mathcal{T}}$ by approximation of its volume averages. Its true volume averages $\average{w}_h=\left[\average{w}\kh\right]_{k=1}^{n_\mathcal{T}}$ are defined on a given finite volume $\mathcal{T}\kh$ by
\begin{align}
\average{w}\kh(t)=\frac{1}{h}\int_{z\khmhalf}^{z\khphalf}w(t,z)\,\mathrm{d}z.
\end{align}
We proceed to semi-discretize the reactor volume equations by identifying midpoint values $w\kh(t)=w(t,z\kh)$ and making extensive use of the midpoint approximation
\begin{align}\label{eq:midpointapprox}
    \average{f(w)}\kh = f(w\kh) +\mathcal{O}(h^2).
\end{align}

\subsubsection{Mass and energy balances}
The average value of the transport equation~\eqref{eq:massbalances} on a finite volume $\mathcal{T}\kh$ is given by
\begin{align}\label{eq:discretizationstart}
    \average{\partial_tc}\kh=-\average{\partial_zN}\kh+\average{R}\kh.
\end{align}
Applying the fundamental theorem of calculus, we arrive at
\begin{align}
    \average{\partial_tc}\kh=-\frac{1}{h}\left(N\khphalf-N\khmhalf\right)+\average{R}\kh,
\end{align}
where $N_h^{k\pm1/2}$ denotes the interface molar fluxes. Applying the midpoint approximation \eqref{eq:midpointapprox}, we see that (up to $\mathcal{O}(h^2)$)
\begin{align}\label{eq:discretizationend}
    \partial_tc\kh=-\frac{1}{h}\left(N\khphalf-N\khmhalf\right)+R\kh,
\end{align}
In a similar fashion, taking the average value of the energy balance \eqref{eq:energybalances} leads to
\begin{align}\label{eq:discretizationenergy}
    \partial_tu\kh=-\frac{1}{h}\left(E\khphalf-E\khmhalf\right)+Q\kh,
\end{align}
where $E_h^{k\pm1/2}$ denotes the interface energy fluxes. The discretized reaction and heat transfer terms, $R\kh$ and $Q\kh$, are defined in Section~\ref{eq:discretetransfer}. When summed over the entire domain, \eqref{eq:discretizationend}-\eqref{eq:discretizationenergy} satisfy discrete analogs of the conservation laws for mass and energy. Moreover, they approximate the corresponding integral balances up to $\mathcal{O}(h)$ in space.

\subsubsection{Upwind transport}
We consider the discretization of the transport models defined in Sections~\ref{sec:masstransport}–\ref{sec:heattransport}. The interfacial fluxes are computed using upwind-based numerical flux approximations determined by the local flow direction, as outlined below. In general, the interface fluxes
\begin{subequations}
\begin{align}
    N\khphalf&=N(w\khphalf,\partial_zw\khphalf),\\
    E\khphalf&=E(w\khphalf,\partial_zw\khphalf),
\end{align}
\end{subequations}
are written in terms of all interface state variables
\begin{align}
    w\khphalf&=[c\khphalf;u\khphalf;T\khphalf;P\khphalf],
\end{align}
and their similarly defined interface gradients $\partial_zw\khphalf$. Since the fluid phase flows in the direction of the negative pressure gradient, we adopt an upwind scheme based on the local pressure difference
\begin{align}
    w\khphalf=\begin{cases}
        w\kh, \quad &\text{if } P\khp \leq P\kh,\\
        w\khp, \quad &\text{if } P\khp > P\kh.
    \end{cases}
\end{align}
In contrast, the interface gradients are resolved directly
\begin{align}\label{eq:gradientresolve}
    \partial_zw\khphalf=\frac{1}{h}\left(w\khp-w\kh\right),
\end{align}
as described in \cite{leveque2002finite}.

\noindent\textbf{Molar fluxes}\\
\noindent The discrete formulation mirrors the decomposition of the molar flux into its advective and diffusive parts
\begin{align}
    N\khphalf = (N^\text{advection})\khphalf+(N^\text{diffusion})\khphalf.
\end{align}
The advective flux is computed using the interface evaluation of the advection model
\begin{align}
    (N^\text{advection})\khphalf &= v\khphalf c\khphalf,
\end{align}
Here, the interface velocity is defined by
\begin{subequations}
\begin{align}\label{eq:discretevelocity}
    v\khphalf &= v(\partial_z P\khphalf,\mu\khphalf,\rho\khphalf),\\
    \mu\khphalf&=\mu(T\khphalf,P\khphalf,c\khphalf),\\
    \rho\khphalf&=\rho(c\khphalf).
\end{align}
\end{subequations}
In a similar fashion, the diffusive flux is computed by evaluating the diffusion model
\begin{align}\label{eq:diff1}
    (N^\text{diffusion})\khphalf &= -D\khphalf \odot\partial_zc\khphalf,
\end{align}
with interface diffusion coefficients
\begin{align}\label{eq:diff2}
    D\khphalf&=D(T\khphalf,P\khphalf,c\khphalf).
\end{align}

\noindent\textbf{Energy flux}\\
\noindent The same procedure applies to the discrete energy flux, which we split into its convective and conductive parts
\begin{align}
    \hspace{-0.5em}E\khphalf = (E^\text{convection})\khphalf+(E^\text{conduction})\khphalf.
\end{align}
The convective part is handled by the thermodynamic model
\begin{multline}
    (E^\text{convection})\khphalf =\\
     H(T\khphalf,P\khphalf,(N^\text{advection})\khphalf),
\end{multline}
possibly with an $\varepsilon$ scaling for fluid-solid volumes. The conductive part is handled analogously
\begin{multline}
    (E^\text{conduction})\khphalf =\\ \phantom{ (E^\text{conduction})\khphalf}-(\bar{H}\khphalf)^T(N^\text{diffusion})\khphalf
    \\-\kappa\khphalf \partial_zT\khphalf,
\end{multline}
with interface partial molar enthalpies and conductivity
\begin{subequations}
\begin{align}
    \bar{H}\khphalf &= \partial_{c}H(T\khphalf,P\khphalf,c\khphalf),\\
    \kappa\khphalf&=\kappa(T\khphalf,P\khphalf,c\khphalf).
\end{align}
\end{subequations}

\subsubsection{Reaction and heat transfer}\label{eq:discretetransfer}
Both the reaction and heat transfer terms are local to the volume. Accordingly, they are directly discretized using the midpoint approximation \eqref{eq:midpointapprox}. The resulting expressions for the discrete production term are summarized by
\begin{subequations}
\begin{align}
    R\kh &= \nu^Tr\kh,\\
    r\kh&=r(T\kh,P\kh,c\kh),
\end{align}
\end{subequations}
and likewise, the discrete heat transfer is defined by
\begin{subequations}
\begin{align}
    Q\kh&=aK\kh\left(T_h^{\prime k}-T\kh\right),\\
    K\kh&=K(T_h^{\prime k},P_h^{\prime k},c_h^{\prime k};T\kh,P\kh,c\kh),
\end{align}
\end{subequations}
with respect to an adjacent reactor volume $V'\in\mathcal{V}$.

\subsubsection{Thermodynamic constraints}\label{eq:discretethermo}
Lastly, we consider the discretization of the thermodynamic constraints proposed by equations \eqref{eq:volumeconstraint}-\eqref{eq:energyconstraint}. Similar to the above, they are discretized by midpoint evaluation
\begin{subequations}
\begin{align}
    V(T\kh,P\kh,c\kh)&=1,\\
    U(T\kh,P\kh,c\kh)&=u\kh.
\end{align}
\end{subequations}

\subsubsection{Initial and boundary conditions}
\noindent\textbf{Initial conditions}\\
\noindent The initial conditions are also discretized by midpoint evaluation, setting
\begin{subequations}\label{eq:discreteIC}
\begin{align}
    c\kh(0)&=c_0(z\kh),\\
    u\kh(0)&=u_0(z\kh),\\
    T\kh(0)&=T_0(z\kh),\\
    P\kh(0)&=P_0(z\kh).
\end{align}
\end{subequations}
\noindent\textbf{Boundary conditions}\\
\noindent  
Given prescribed flow-driven inflow BCs, the inlet conditions are simply implemented by the expressions in \eqref{eq:flowinlet}
\begin{subequations}
\begin{align}
    N_h^{1/2}(t)=\frac{f_\text{in}(t)}{S_\text{fluid}},\\
    E_h^{1/2}(t)=\frac{h_\text{in}(t)}{S}.
\end{align}
\end{subequations}
In contrast, implementing pressure-driven inflow BCs requires careful consideration. They are given by
\begin{subequations}
\begin{align}
    N_h^{1/2}(t) &=v_h^{1/2}(t) c_\text{in}(t),\\
    E_h^{1/2}(t) &=H(T_\text{in}(t),P_\text{in}(t),N_h^{1/2}(t)),
\end{align}
\end{subequations}
possibly with an $\varepsilon$ scaling of the enthalpy flux for fluid-solid volumes. Towards their implementation, we define the inlet velocity as in \eqref{eq:flowproperties}
\begin{subequations}\label{eq:discreteflowproperties}
\begin{align}
    v_h^{1/2}(t) &= v(\partial_z P_h^{1/2}(t),\mu_\text{in}(t),\rho_\text{in}(t)),\\
    \mu_\text{in}(t)&=\mu(T_\text{in}(t),P_\text{in}(t),c_\text{in}(t)),\\
    \rho_\text{in}(t)&=\rho(c_\text{in}(t)),
\end{align}
\end{subequations}
given the inlet pressure gradient computed from the prescribed inlet pressure
\begin{align}
    \partial_z P_h^{1/2}(t)&=\frac{P_h^1(t)-P_\text{in}(t)}{h/2}.
\end{align}
The free outflow BCs are implemented in a similar fashion using the upwind flux approximation
\begin{subequations}
\begin{align}
    N_h^{n_\mathcal{T}+1/2}(t) &=v_h^{n_\mathcal{T}+1/2}(t) c_h^{n_\mathcal{T}}(t),\\
    E_h^{n_\mathcal{T}+1/2}(t) &=H(T_h^{n_\mathcal{T}}(t),P_h^{n_\mathcal{T}}(t),N_h^{n_\mathcal{T}+1/2}(t)),
\end{align}
\end{subequations}
with outlet velocity
\begin{subequations}
\begin{align}
    v_h^{n_\mathcal{T}+1/2}(t) &= v(\partial_z P_h^{n_\mathcal{T}+1/2}(t),\mu_h^{n_\mathcal{T}}(t),\rho_h^{n_\mathcal{T}}(t)),\\
    \mu_h^{n_\mathcal{T}}(t)&=\mu(T_h^{n_\mathcal{T}}(t),P_h^{n_\mathcal{T}}(t),c_h^{n_\mathcal{T}}(t)),\\
    \rho_h^{n_\mathcal{T}}(t)&=\rho(c_h^{n_\mathcal{T}}(t)),
\end{align}
\end{subequations}
defined by the outlet pressure gradient given by the prescribed outlet pressure
\begin{align}
    \partial_z P_h^{n_\mathcal{T}+1/2}(t)=\frac{P_\text{out}(t)-P_h^{n_\mathcal{T}}(t)}{h/2}.
\end{align}

\subsection{Semi-discrete model}
We denote the discretized differential and algebraic state variables for a volume $V \in \mathcal{V}$ by vectors $x_h~=[c_h;u_h]$ and $y_h~=[T_h;P_h]$, and combine all volume variables into differential and algebraic state vectors $\boldsymbol{x}$ and $\boldsymbol{y}$, respectively. Accordingly, the spatial discretization can be summarized by a semi-discrete system of semi-explicit index-1 differential-algebraic equations (DAEs) \citep{hairer1996solving}
\begin{subequations} \label{eq:semidiscrete}
\begin{align}
    \dot{\boldsymbol{x}}&=f(\boldsymbol{x},\boldsymbol{y},p), \quad &\boldsymbol{x}(0)&=\boldsymbol{x}_0,\\
    0&=g(\boldsymbol{x},\boldsymbol{y},p), \quad &\boldsymbol{y}(0)&=\boldsymbol{y}_0,
\end{align}
\end{subequations}
with initial conditions $(\boldsymbol{x}_0,\boldsymbol{y}_0)$ defined by their discretization \eqref{eq:discreteIC}, and where $p$ denotes a vector of parameters, e.g., boundary conditions, submodel parameters, and possible model inputs or disturbances.

\subsubsection{Model derivatives}
To support simulation, we require derivatives of the semi-discrete model functions in \eqref{eq:semidiscrete}. Derivatives of the thermodynamic model were discussed in Section~\ref{subsec:thermo}, but similar strategies apply to the remaining constitutive functions. Indeed, if analytical differentiation is impractical, automatic differentiation (AD) \citep{griewank2008evaluating} can be used via open-source software libraries \citep{revels2016forwardmode,jax2018github}. The spatial discretization in Section~\ref{subsec:spatial} also admits analytical derivatives, but its sparsity structure makes it well suited for efficient sparse AD \citep{hill2025sparser}.

\subsection{Steady-state simulation}
Steady-state reactor profiles can be obtained by setting $\dot{\boldsymbol{x}}=0$  in \eqref{eq:semidiscrete} and solving the resulting system
\begin{align}
    0=F(\boldsymbol{w},p),
\end{align}
with $F=[f;g]$, for the combined vector $\boldsymbol{w}=[\boldsymbol{x};\boldsymbol{y}]$, at fixed conditions $p$. To this end, we employ Newton's method
\begin{subequations}
\begin{align}
    \partial_{\boldsymbol{w}} F(\boldsymbol{w}^\ell,p)\Delta
    \boldsymbol{w}^{\ell}&=-F(\boldsymbol{w}^\ell,p),\label{eq:newtonsys}\\
    \boldsymbol{w}^{\ell+1}&=\boldsymbol{w}^{\ell}+\alpha^{\ell}\Delta\boldsymbol{w}^\ell,
\end{align}
\end{subequations}
where the step size $\alpha^\ell$ is chosen by line search globalization \citep{nocedal2006numerical}. Starting from an initial guess $\boldsymbol{w}^0$, we iterate over $\ell = 0, 1, 2, \dots$, until convergence. The residual Jacobian $\partial_{\boldsymbol{w}} F$ is computed by the model derivatives, and the Newton system \eqref{eq:newtonsys} may be efficiently solved by a sparse direct linear solver \citep{davis2004algorithm,davis2010algorithm}.

\subsubsection{Initial guess}\label{subsec:initialguess}
Despite incorporating a globalization procedure, it is important to supply a sufficiently accurate initial guess to ensure fast convergence speed, i.e., reasonable simulation times. In this work, we initialize reactor volumes by prescribing constant temperature profiles (assuming no reaction), linear pressure profiles from a prescribed pressure drop, and component concentrations that satisfy the volume constraint. The internal energy profile is then computed explicitly to enforce the energy constraint, ensuring a thermodynamically consistent starting point. When multiple volumes are connected and an overall pressure drop is prescribed, their respective pressure profiles are initialized by distributing the pressure drop according to the coupling condition in \eqref{eq:convcoupling}, thereby determining consistent volume-specific pressure gradients.

\subsubsection{Numerical continuation}
It may not be possible to generate good initial guesses of steady states under all considered process conditions ${p \in [p_\text{min},p_\text{max}]}$, especially when multiple steady states arise. In particular for parametric sensitivity analysis, standard parameter sweeps near turning points may fail to trace the full steady-state solution manifold. To overcome this, we employ pseudo-arclength continuation (PLAC) \citep{allgower1990numerical} to trace the full solution curve $(\boldsymbol{w}(s), p(s))$ for a single parameter $p$ as a function of the continuation parameter $s$. Given an initial steady-state solution $F(\boldsymbol{w}^0,p^0)=0$, PLAC solves the sequence of augmented systems
\begin{subequations}
\begin{align}
    \hspace{-0.5em}0&=F(\boldsymbol{w}^{m+1},p^{m+1}),\\
    \hspace{-0.5em}0&=(\boldsymbol{w}^{m+1}-\boldsymbol{w}^m)^T\xi_{\boldsymbol{w}}^{m} + (p^{m+1} - p^m)\xi_p^m- \Delta s^{m}
\end{align}   
\end{subequations}
using Newton’s method at each continuation step $m=0,1,2,\dots$. Here, $\Delta s^m$ denotes an adaptively chosen step size, and $\xi^m=[\xi_{\boldsymbol{w}}^{m};\xi_p^m]$ denotes the normalized approximate tangent vector obtained by previous secant information
\begin{align}
    (\boldsymbol{w}^{m},p^{m})=(\boldsymbol{w}^{m-1},p^{m-1})+\Delta s^{m-1}\xi^m.
\end{align}
The first secant approximation is obtained from two solutions that are obtained via natural continuation.

\subsection{Dynamic simulation}
Dynamic reactor trajectories are obtained by solving the semi-discrete DAE system \eqref{eq:semidiscrete}. For reasons that will soon become apparent, we use the combined vector notation $\boldsymbol{w}$ and solve the DAE system in mass matrix form
\begin{align}\label{eq:massmatrix}
    M\dot{\boldsymbol{w}}=F(\boldsymbol{w},p), \quad \boldsymbol{w}(0)=\boldsymbol{w}_0.
\end{align}
Here, the mass matrix $M$ indicates the presence or absence of time derivatives, with which the initial condition $\boldsymbol{w}_0$ is assumed to be consistent, i.e., $F(\boldsymbol{w}_0,p)\in \text{range}\,M$.

\subsubsection{Mass matrices}
The dynamic differential mass and energy balances in \eqref{eq:semidiscrete} are represented by the mass matrix
\begin{align}
    M=\left[\begin{array}{cc}
         I_{n_{\boldsymbol{x}}\times n_{\boldsymbol{x}}}& 0_{n_{\boldsymbol{x}}\times n_{\boldsymbol{y}}} \\
         0_{n_{\boldsymbol{y}}\times n_{\boldsymbol{x}}}& 0_{n_{\boldsymbol{y}}\times n_{\boldsymbol{y}}}
    \end{array}\right],
\end{align}
and consistency is assured by assuming $g(\boldsymbol{x}_0,\boldsymbol{y}_0,p)=0$. Note that the initial guesses described in Section \ref{subsec:initialguess}, as well as any computed steady state, will always be consistent. Following the derivation in Appendix \ref{app:energybalance}, pseudo-steady-state assumptions give rise to different mass matrices, which motivates the mass matrix formulation due to its flexibility in selecting differential state variables without specifying new model functions.

\subsubsection{Diagonally-implicit Runge-Kutta methods}
Towards the numerical solution of \eqref{eq:massmatrix}, we employ the one-step family of diagonally implicit Runge-Kutta (DIRK) methods \citep{hairer1996solving}. Generally, an $s$-stage DIRK method can be formulated by stepping from $(t_n,\boldsymbol{w}_n)$ to $(t_{n+1},\boldsymbol{w}_{n+1})$, with adaptive step size $h_n=t_{n+1}-t_n$, by solving for the stage states $\boldsymbol{W}\!_{n}=[\boldsymbol{W}\!_{n,i}]_{i=1}^{s}$ defined by the implicit stage equations
\begin{align}\label{eq:stageequation}
    M\boldsymbol{W}\!_{n,i}&=M\boldsymbol{w}_n + h_n\sum_{j=1}^s a_{i,j} F(\boldsymbol{W}\!_{n,j},p),
\end{align}
at stage times $T_{n,i}=t_n+c_ih$, for $i=1,\dots,s$. The stage states are combined in a quadrature rule to obtain the next step
\begin{align}\label{eq:quadraturerule}
    M\boldsymbol{w}_{n+1}=M\boldsymbol{w}_{n}+ h_n\sum_{i=1}^s b_i F(\boldsymbol{W}\!_{n,i},p).
\end{align}
The method is characterized by its Butcher tableau coefficients $a=[a_{i,j}]_{i,j=1}^{s}$, $b=[b_i]_{j=1}^s$, and $c=[c_i]_{i=1}^s$, which are chosen to ensure high-order accuracy, stability, and stiff accuracy \citep{kennedy2016diagonally}. In DIRK methods, the lower triangular structure of the coefficient matrix, where $a_{i,j}=0$ for $i<j$, enables sequential solution of the stage equations. We specifically employ singly-diagonally implicit Runge-Kutta (SDIRK) methods, with constant diagonal entries $a_{i,i}=\gamma$, and use explicit first-stage SDIRK (ESDIRK) methods with $a_{1,1}=0$. Using Newton's method to solve \eqref{eq:stageequation} with sparse linear algebra, under these considerations, results in a highly efficient numerical method. Additionally, we adopt stiffly accurate schemes where $b_i=a_{s,i}$, so the final stage directly yields the step solution, rendering the quadrature update \eqref{eq:quadraturerule} redundant. This combination of method properties (stiffly accurate ESDIRK) allows us to efficiently handle index-1 DAE systems.
\subsection{Further applications}
By leveraging the structured DAE formulation, the presented simulation framework readily extends to tasks such as state estimation, and steady-state or dynamic optimization, supporting advanced process design and flexible control strategies, including model predictive control \citep{rawlings2022model, oliveiracabral2024learningbased, kong2024nonlinear}.

%% file: tex/4casestudy.tex
\section{Ammonia synthesis case study}\label{sec:casestudy}
\begin{table}[tb]
\caption{Volume dimensions for both the AFBR and IDCR reactor units.}
\centering
\begin{tabular}{lll|r|rr}
         &               &       & \multicolumn{1}{c|}{AFBR} & \multicolumn{2}{c}{IDCR}                          \\ \cline{4-6} 
Name     & Symbol        & Unit  & \multicolumn{1}{c|}{FBR}  & \multicolumn{1}{c}{FBR} & \multicolumn{1}{c}{LCT} \\ \hline
Length   & $L$           & m     & 2                         & 6                       & 6                       \\
Volume   & $V$           & m$^3$ & 2                         & 3                       & 0.5                     \\
Porosity & $\varepsilon$ & --    & 0.33                      & 0.18                    & --                     
\end{tabular}
\label{tab:reactordimensions}
\end{table}
\begin{table}[tb]
\caption{Parameters for the reactor models' constitutive functions.}
\centering
\resizebox{\columnwidth}{!}{%
\begin{tabular}{lll|r}
Model                          & Symbol               & Unit                              & Value                \\ \hline
Solid                          & $\rho_\text{solid}$  & kg/m$^3$                          & 3284                 \\
                               & $c_{p,\text{solid}}$ & J/(kg$\cdot$K)                    & 1100                 \\ \hline
Kinetic                        & $\eta$               & --                                & 4.75                 \\
                               & $\beta$              & --                                & 0.5                  \\
                               & $A_\rightarrow$      & mol/(s$\cdot$m$^3$-solid) & 4\,972               \\
                               & $A_\leftarrow$       & mol/(s$\cdot$m$^3$-solid) & 7.14$\cdot$10$^{15}$ \\
                               & $E_\rightarrow$      & J/mol                             & 87\,090              \\
                               & $E_\leftarrow$       & J/mol                             & 198\,464             \\
                               & $R_\text{ideal}$     & J/(mol$\cdot$K)                   & 8.314                \\ \hline
Advection                      & $\mu$                & Pa$\cdot$s                        & 3.08$\cdot$10$^{-5}$ \\
                               & $d_p$                & m                                 & 8$\cdot$10$^{-3}$    \\
                               & $d_t$                & m                                 & 13.3$\cdot$10$^{-3}$ \\
                               & $f_{DW}$             & --                                & 2$\cdot$10$^{-2}$    \\ \hline
Diffusion                      & $D$                  & m$^2$/s                           & 1$\cdot$10$^{-5}$    \\ \hline
Conduction                     & $\kappa$             & W/(m$\cdot$K)                     & 50                 \\ \hline
Heat transfer & $\upsilon$           & W/(m$^2\cdot$K)                   & 300                  \\
              & $A$                  & m$^2$                             & 150                 
\end{tabular}%
}
\label{tab:reactorparameters}
\end{table}

\begin{table}[tb]
\centering
\caption{Nominal operating conditions given pressure-based BCs.}
\begin{tabular}{l|cccc|cc}
Symbol & $x_{\text{in},\text{N}_2}$ & $x_{\text{in},\text{H}_2}$ & $x_{\text{in},\text{NH}_3}$ & $x_{\text{in},\text{Ar}}$ & $P_\text{in}$ & $P_\text{out}$ \\ 
Unit   & \%        & \%         & \%          & \%        & bar     & bar      \\ \hline
Value  & 21.5           & 64.5           & 10              & 4             & 200           & 199           
\end{tabular}%
\label{tab:reactorconditions}
\end{table}

In this section, we demonstrate the capabilities of the proposed modeling and simulation framework by applying it to the synthesis of ammonia. We consider the Haber-Bosch (HB) process, which synthesizes ammonia via the exothermic reaction
\begin{align}\label{eq:ammoniasynthesis}
    \text{N}_2+3\text{H}_2 \rightleftharpoons 2\text{NH}_\text{3}.
\end{align}
The reaction is carried out over iron-based catalysts at elevated temperatures $T\in[500,800]$ [K] and high pressures $P\in[150,300]$ [bar] to balance the kinetic and thermodynamic reaction limitations. Due to the equilibrium limitations, single-pass conversions are typically limited to 10–30\%, necessitating product separation via condensation and extensive recycle of unreacted gases to achieve overall conversions exceeding 95\% \citep{nielsen1995ammonia}. In the context of P2A production, the HB process is powered by renewable electricity sources such as wind or solar. Hydrogen is supplied via water electrolysis, while nitrogen is obtained through air separation, preferably by pressure swing adsorption (PSA) or membrane-based technologies, which offer great operational flexibility \citep{palys2022powertox}.

To investigate the behavior of the reactor units under varying model assumptions and operating conditions relevant to green ammonia synthesis, we develop models for two representative reactor units. We model both reactors using the framework established in Section~\ref{sec:modeling} with the reactor dimensions summarized in Table~\ref{tab:reactordimensions}. First, we consider an adiabatic fixed-bed reactor (AFBR), typical of individual beds in multi-bed reactor configurations (e.g., AQCRs or AICRs), and second, an isothermal direct-cooled reactor (IDCR) incorporating (autothermal) counter-current feed gas cooling. For both configurations, we investigate the impact of employing real versus ideal fluid thermodynamics, dispersive versus non-dispersive transport, and fully dynamic versus pseudo-steady-state assumptions.

\subsection{Reactor models}
The AFBR is modeled by the FBR equations presented in Section \ref{subsec:FBR}, assuming adiabatic operation $(Q=0)$. On the other hand, the IDCR is modeled using the DCR equations from Section \ref{subsec:DCR}, coupled via the volume-coupling condition \eqref{eq:couplingconditions}, which connects the outlet of the LCT volume to the inlet of the FBR volume. The common model components listed in Table~\ref{tab:modelcomponents} are defined using the parameters provided in Table~\ref{tab:reactorparameters} and specified below.

\subsubsection{Thermodynamic model}
\noindent\textbf{Fluid model}\\
\noindent We employ the polynomial ideal-gas heat capacity correlations provided in the appendix of \cite{poling2000properties}, and consider both ideal-gas and real-gas behavior using the cubic EOS by Soave–Redlich–Kwong (SRK) \citep{soave1972equilibrium}. 

\noindent\textbf{Solid model}\\
The solid-phase internal energy density is computed assuming negligible pressure-work and constant solid properties. Given the reference state $u_\text{solid}(T_\text{ref})=0$ at $T_\text{ref}=0$, we let
\begin{align}\label{eq:solidenergy}
    u_\text{solid}(T)=\rho_\text{solid}c_{p,\text{solid}}T,
\end{align}
where $c_{p,\text{solid}}$ [J/(kg$\cdot$K)] and $\rho_\text{solid}$ [kg/m$^3$] denote the solid's specific heat capacity and density, respectively. They are specified as in \cite{rosbo2024comparison}.

\subsubsection{Reaction model}
\noindent\textbf{Stoichiometric model}\\
\noindent
We consider the components $\mathcal{C}=\{\text{N}_2,\text{H}_2,\text{NH}_3,\text{Ar}\}$, in which we include argon as a trace inert from upstream air separation. The ammonia reaction \eqref{eq:ammoniasynthesis} is summarized by the stoichiometric matrix
\begin{align}
    \nu=[-1, -3, \phantom{+}2, \phantom{+}0].
\end{align}
\noindent\textbf{Kinetic model}\\
\noindent
The kinetic rate is described by the Temkin-Pyzhev rate expression \citep{froment2010chemical}
\begin{align}\label{eq:temkin}
    r=k_\rightarrow P_{\text{N}_2}\left(\frac{P^{3}_{\text{H}_2}}{P^2_{\text{NH}_3}}\right)^\beta-k_\leftarrow \left(\frac{P^2_{\text{NH}_3}}{P^3_{\text{H}_2}}\right)^\beta.
\end{align}
Here $\beta$ [-] is a fixed exponent, and $[P_\alpha]_{\alpha\in\mathcal{C}}$ [bar] are partial pressures defined by Dalton's law, $P_\alpha=Px_\alpha$, where $x=[x_\alpha]_{\alpha\in\mathcal{C}}$ [-] denotes the component mole fractions which can be computed by $x_\alpha=c_\alpha/\sum_\alpha c_\alpha$. The forward and backward rate constants are described by Arrhenius expressions
\begin{subequations}
\begin{align}
    k_\rightarrow(T)&= A_\rightarrow \exp\left(\frac{E_\rightarrow}{R_\text{ideal}T}\right),\\
    k_\leftarrow(T)&= A_\leftarrow \exp\left(\frac{E_\leftarrow}{R_\text{ideal}T}\right),
\end{align}
\end{subequations}
 where $A_\rightarrow,A_\leftarrow$ [mol/(s$\cdot$m$^3$-solid)] denotes the Arrhenius factors, and $E_\rightarrow,E_\leftarrow$ [J/mol] denotes the activation energies for both reaction directions. Additionally, $R_\text{ideal}$ [J/(mol$\cdot$K)] denotes the ideal gas constant. The rate expression \eqref{eq:temkin} is scaled to include the pseudo-homogeneous scaling factors from \eqref{eq:pseudohomreaction}. However, note that the effectiveness factor is considered a fitted activity factor that corrects the kinetic rates to account for the effectiveness of modern-day catalysts \citep{morud1998analysis,rosbo2024comparison}.

\subsubsection{Mass transport model}
\noindent\textbf{Advection model}\\
\noindent
We model the pressure drop inside the FBR volumes using the Ergun equation \citep{froment2010chemical}
\begin{align}\label{eq:erguntext}
    -\partial_zP&=\frac{150\mu(1-\varepsilon)^2}{d_\text{p}^2\varepsilon^2}v\\
    &+\frac{1.75\rho_\text{fluid}(1-\varepsilon)}{d_\text{p}\varepsilon}v|v|,\nonumber
\end{align}
where $d_\text{p}$ [m] denotes the catalyst particle diameter. The fluid viscosity $\mu$ is set constant to that of air at 600 [K], and the fluid density $\rho_\text{fluid}=\rho_\text{fluid}(c)$ is calculated using~\eqref{eq:fluiddensity}. Moreover, we model the pressure drop inside the LCT volume using the Darcy-Weisbach equation
\begin{align}\label{eq:darcyweisbachtext}
    -\partial_zP=\frac{f_\text{DW}\rho_\text{fluid}}{2d_\text{t}}v|v|,
\end{align}
where $d_\text{t}$ [m] denotes the individual cooling tube diameter, and $f_\text{DW}$~[-] denotes the Darcy-Weisbach friction factor chosen in the turbulent flow regime from a Moody diagram \citep{bird2006transport}. We invert the expressions \eqref{eq:erguntext}-\eqref{eq:darcyweisbachtext} to explicit functions for the molar average velocities following~\cite{martinsen2023modeling}.

\noindent\textbf{Diffusion model}\\
\noindent
A common effective diffusion coefficient is chosen as a constant average value following the correlations presented in \cite{martinsen2023modeling}.

\subsubsection{Energy transport model}
\noindent\textbf{Conduction model}\\
\noindent
The solid thermal conductivity is modeled as a constant representative average value for iron in the temperature range $T\in[600,800]$ [K] \citep{engineeringtoolbox2025thermal}. In comparison, the gas thermal conductivity is negligible.

\subsubsection{Heat transfer model}
Relevant only for the IDCR, the overall heat transfer coefficient is approximated as a constant average value for high-pressure gases inside and outside tubes \citep{engineeringtoolbox2025heat}.

\subsubsection{Initial and boundary conditions}
We prescribe pressure-driven BCs with a total pressure drop of $\Delta P=1$ [bar]. The nominal operating conditions are specified in Table \ref{tab:reactorconditions}.

\subsection{Steady-state and dynamic simulation}
We discretize each reactor volume into $n_\mathcal{T}=100$ finite volumes. Steady-state and dynamic simulations are conducted at varying inlet temperatures using the methodology presented in Section \ref{sec:simulation}. Initial guesses for steady-state solutions are found and subsequently refined using a line-search Newton method, before being swept across a range of inlet temperatures by PLAC. The resulting set of steady states serve as consistent initial conditions for dynamic simulation, which is performed using an ESDIRK method. All numerical methods are supplied with equal absolute and relative tolerances. We set the tolerance level $\text{tol}=10^{-4}$ for steady-state simulations and $\text{tol}=10^{-5}$ for dynamic simulations.

\begin{figure*}[tb]
    \centering
    \begin{subfigure}{0.325\textwidth}
        \includegraphics[width=\textwidth]{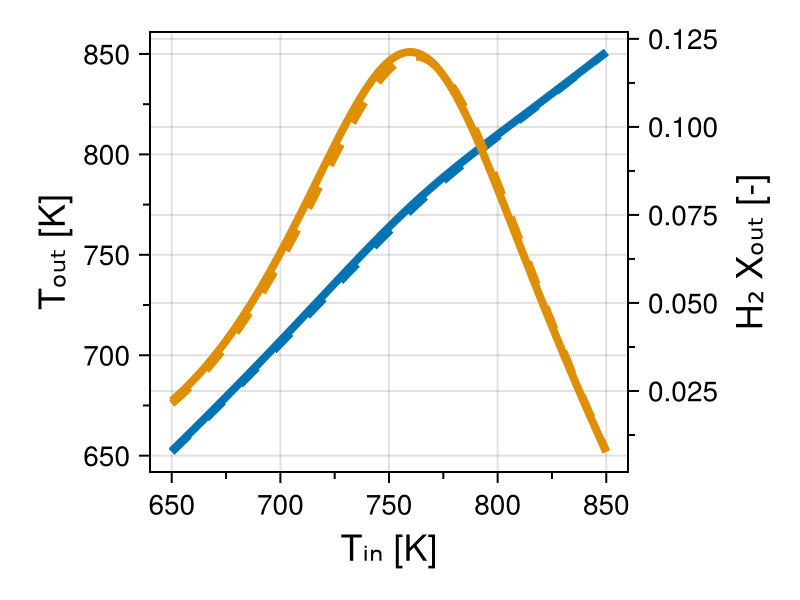}
        \caption{}
        \label{subfig:AFBRnominalcurves}
    \end{subfigure}
    \hfill
    \begin{subfigure}{0.325\textwidth}
        \includegraphics[width=\textwidth]{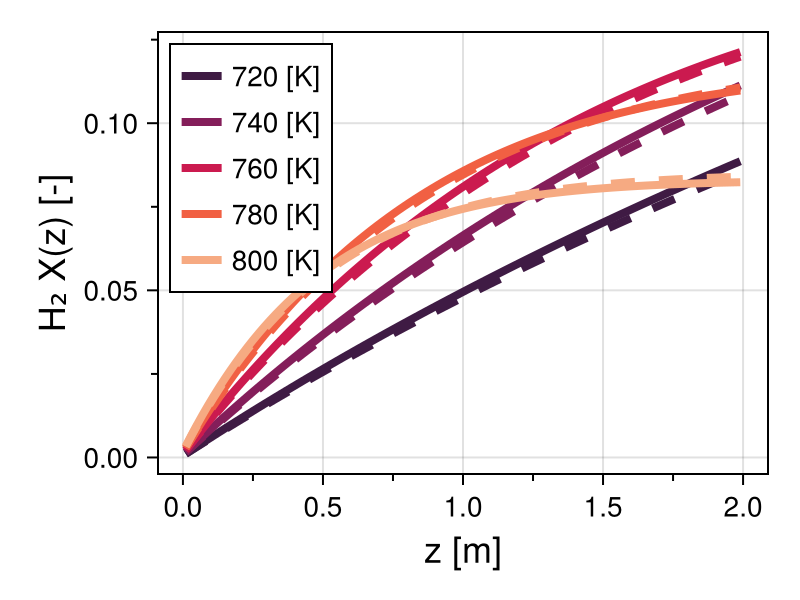}
        \caption{}
        \label{subfig:AFBRnominalconversionprofiles}
    \end{subfigure}
    \hfill
    \begin{subfigure}{0.325\textwidth}
        \includegraphics[width=\textwidth]{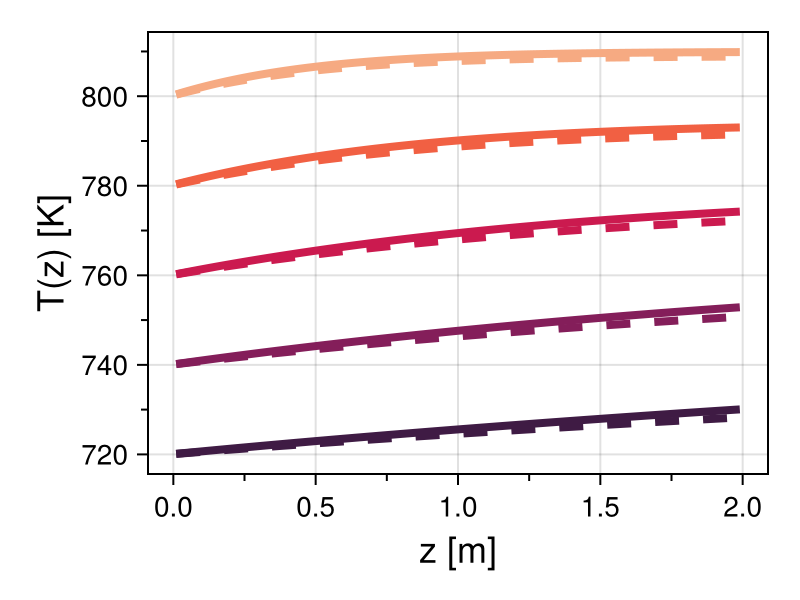}
        \caption{}
        \label{subfig:AFBRnominaltemperatureprofiles}
    \end{subfigure}
    \caption{AFBR steady states at nominal operating conditions. \textbf{(a):} Outlet H$_2$ conversion (yellow) and outlet temperature~(blue). \textbf{(b):} H$_2$ conversion profiles. \textbf{(c):} Temperature profiles. \textbf{(a)--(c):} Real-fluid properties (solid) and ideal-fluid properties (dashed).}
    \label{fig:AFBRnominal}
\end{figure*}

\begin{figure}[tb]
    \centering
    \begin{subfigure}{\linewidth}
        \centering
        \includegraphics[width=0.75\linewidth]{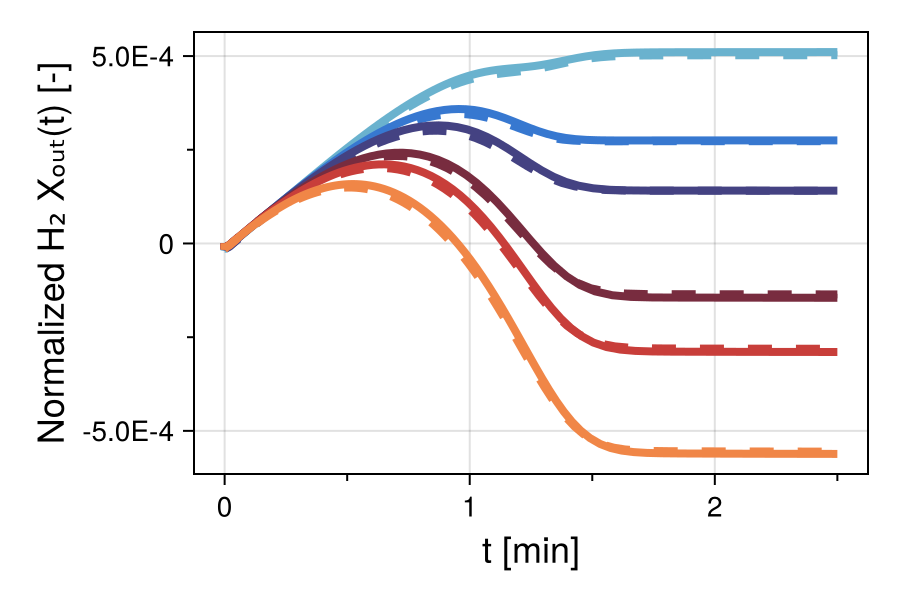}
        \caption{}
        \label{subfig:AFBRstepsconversion}
    \end{subfigure}
    \begin{subfigure}{\linewidth}
        \centering
        \includegraphics[width=0.75\linewidth]{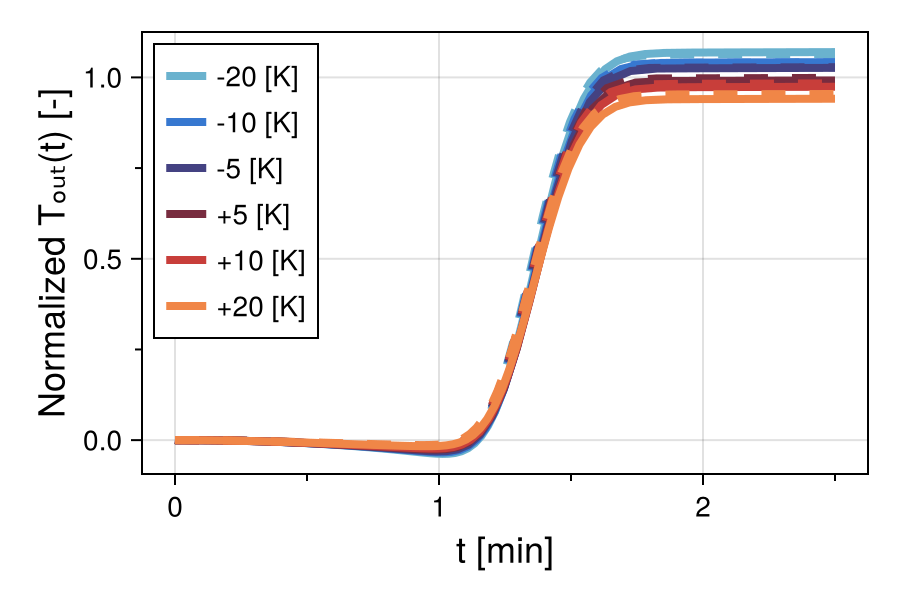}
        \caption{}
        \label{subfig:AFBRstepstemperature}
    \end{subfigure}
    \caption{AFBR step responses for varying steps in the inlet temperature at optimal operating conditions. \textbf{(a):} Normalized H$_2$ conversion. \textbf{(b):} Outlet temperature. \textbf{(a)--(b):} Real-fluid properties (solid) and ideal-fluid properties (dashed).}
    \label{fig:AFBRsteps}
\end{figure}

\subsubsection{Real and ideal fluid thermodynamics}
We investigate the impact of different thermodynamic models on simulated outlet conditions and reactor profiles by comparing ideal fluid properties with real fluid properties described by the cubic SRK EOS.

\noindent\textbf{AFBR}\\
\noindent 
We first perform steady-state simulations of the AFBR under nominal operating conditions with varying inlet temperatures $T_\text{in}~\in~[650,850]$ for both the real and ideal thermodynamic models. The results are seen in Figure~\ref{fig:AFBRnominal}. In Figure~\ref{subfig:AFBRnominalcurves} we see the resulting steady-state outlet temperatures $T_\text{out}$ and outlet H$_2$ conversion $X_{\text{out}}$ plotted against $T_\text{in}$. From these, it is qualitatively hard to discern the impact of the thermodynamic model. Quantitatively, the real model predicts an optimal inlet temperature of roughly $T^*_\text{in}=760$ [K] which is $2$ [K] lower than the optimal inlet temperature for the real model at $T^*_\text{in}=762$ [K]. The corresponding outlet H$_2$ conversion are $X_{\text{out}}^*=12.1\%$ for the real model, and only $0.1\%$ lower for the ideal model. A similar picture is seen by the the steady-state profiles shown in Figures~\ref{subfig:AFBRnominalconversionprofiles}-\ref{subfig:AFBRnominaltemperatureprofiles}. They show that the real model consistently predicts higher reactor temperature profiles, with the largest difference at the outlet. The maximum difference amounts to roughly $\Delta T_\text{out}=2$ [K] at $T_\text{in}=747$ [K], before the predicted $T^*_\text{in}$ for both models. A similar situations holds for the outlet fractions, but the difference between thermodynamic models is negligible.

The steady-state profiles in Figure~\ref{subfig:AFBRnominalconversionprofiles} illustrate how the kinetic and thermodynamic equilibria limit the ammonia reaction. Indeed, below the optimal inlet temperature, the reaction proceeds at a slower, but steady, kinetically-limited rate. However, above the optimal inlet temperature, the reaction initially proceeds rapidly before it is hindered by the thermodynamic equilibrium. This behavior is also seen from the step responses in Figure~\ref{fig:AFBRsteps}. Figure~\ref{subfig:AFBRstepsconversion} shows that increasing the inlet temperature will temporarily increase the conversion, but will ultimately lead to a steady-state with lower conversion. As a result, we see that the outlet temperature response in Figure~\ref{subfig:AFBRstepstemperature} is delayed by approximately 1 [min]. Both models present similar transient behavior.

\begin{figure}[tb]
    \centering
    \includegraphics[width=0.75\linewidth]{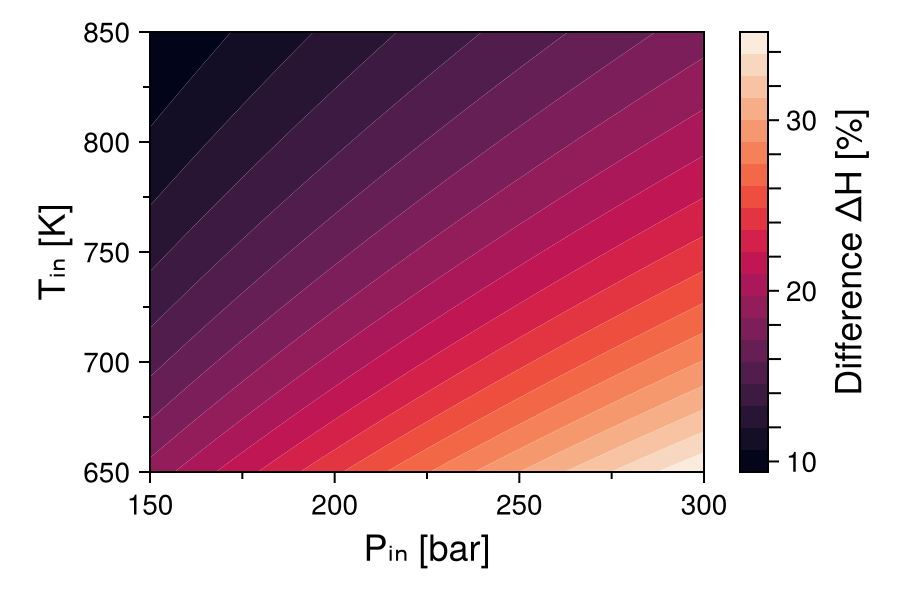}
    \caption{Difference in the predicted heat of reaction $\Delta H$ of the real model against the ideal thermodynamic model.}
    \label{fig:AFBRdeltaH}
\end{figure}

\begin{figure*}[tb]
    \centering
    \begin{subfigure}{0.325\textwidth}
        \includegraphics[width=\textwidth]{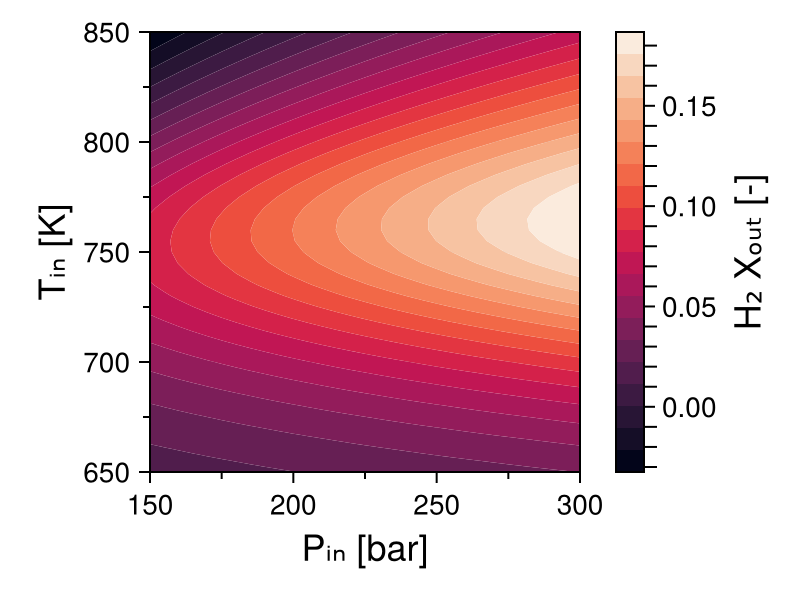}
        \caption{}
        \label{subfig:AFBRsurfacesconversion}
    \end{subfigure}
    \hfill
    \begin{subfigure}{0.325\textwidth}
        \includegraphics[width=\textwidth]{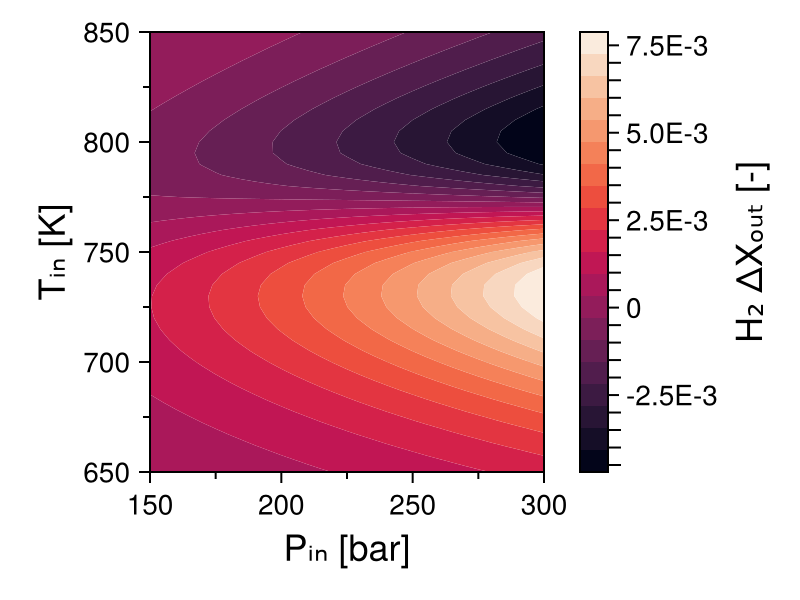}
        \caption{}
        \label{subfig:AFBRsurfacesdeltaconversion}
    \end{subfigure}
    \hfill
    \begin{subfigure}{0.325\textwidth}
        \includegraphics[width=\textwidth]{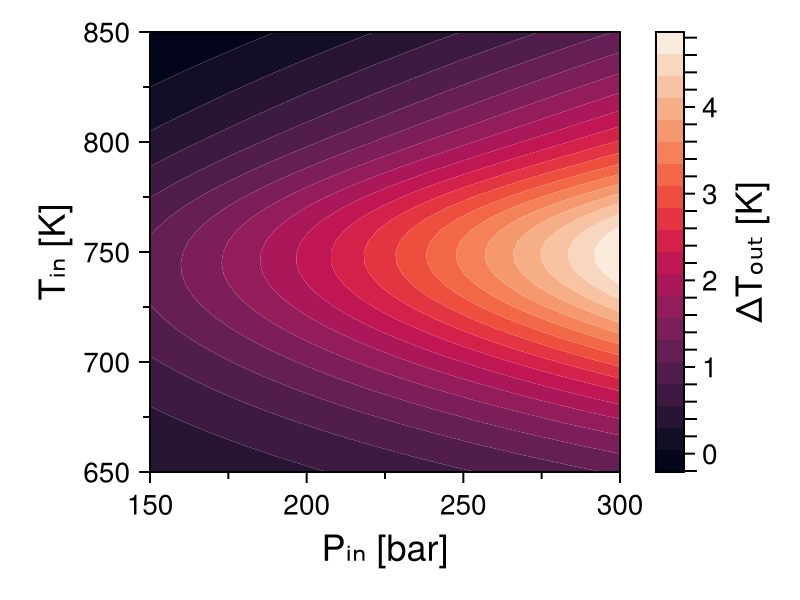}
        \caption{}
        \label{subfig:AFBRsurfacesdeltatemperature}
    \end{subfigure}
    \caption{AFBR steady-state outlets at varying operating conditions. \textbf{(a):} Real model outlet H$_2$ conversion. \textbf{(b):} Difference in outlet H$_2$ conversion of the real model against the ideal thermodynamic model. \textbf{(c):} Difference in outlet temperature.}
    \label{fig:AFBRsurfaces}
\end{figure*}

The operating pressure $P_\text{in} = 200$ [bar] contributes significantly to the thermodynamic equilibrium. To better understand the differences between the thermodynamic models at varying pressures $P_\text{in} \in [150, 300]$, we begin by comparing their predicted heat of reaction $\Delta H$ (as defined in Appendix~\ref{app:energybalance}). The results are shown in Figure~\ref{fig:AFBRdeltaH}. Compared to the ideal model's pressure-independent predictions, the real model predicts stronger exothermic behavior at higher pressures. Accordingly, the largest difference is found at lower temperatures. At the previous optimal operating conditions of $P=200$ [bar] and $T=760$ [K], the real model predicts a heat of reaction that is 16\% more exothermic. The effect of elevated pressures and model difference is demonstrated in Figure~\ref{fig:AFBRsurfaces}. As expected, both models predict higher conversions at elevated pressures. Consistent with our previous observation, the real-fluid model reaches its optimum at slightly lower inlet temperatures. It also predicts higher outlet temperatures, which can be attributed to its more exothermic heat of reaction. Although the difference in the heat of reaction is more pronounced at lower temperatures, its effect is negligible due to kinetic limitations imposed by the large activation energy. For comparison, we repeated the simulations using the cubic Peng–Robinson EOS~\citep{peng1976new}. Relative to the SRK EOS, the predicted heat of reaction $\Delta H$ deviates by at most 1.5\% across similar conditions, explaining the negligible differences observed in their simulation results.

\begin{figure*}[tb]
    \centering
    \begin{subfigure}{0.325\textwidth}
        \includegraphics[width=\textwidth]{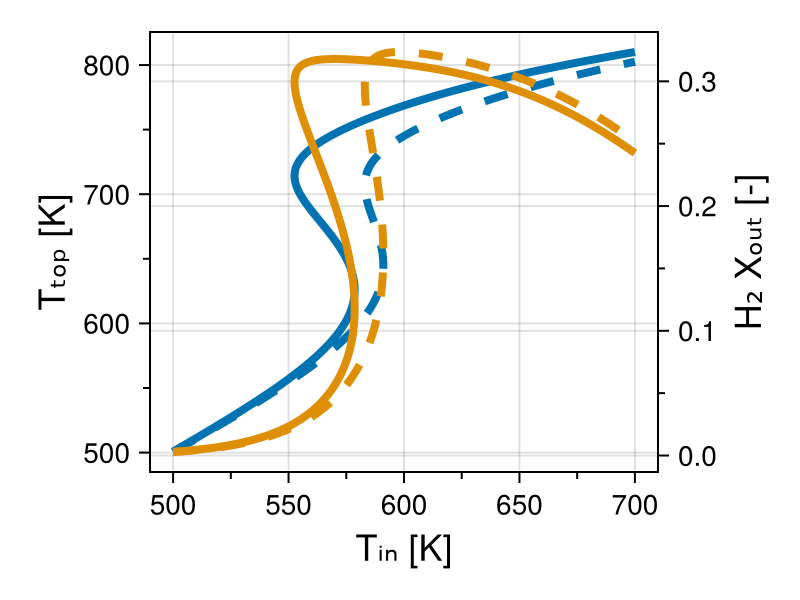}
        \caption{}
        \label{subfig:IDCRnominalcurves}
    \end{subfigure}
    \hfill
    \begin{subfigure}{0.325\textwidth}
        \includegraphics[width=\textwidth]{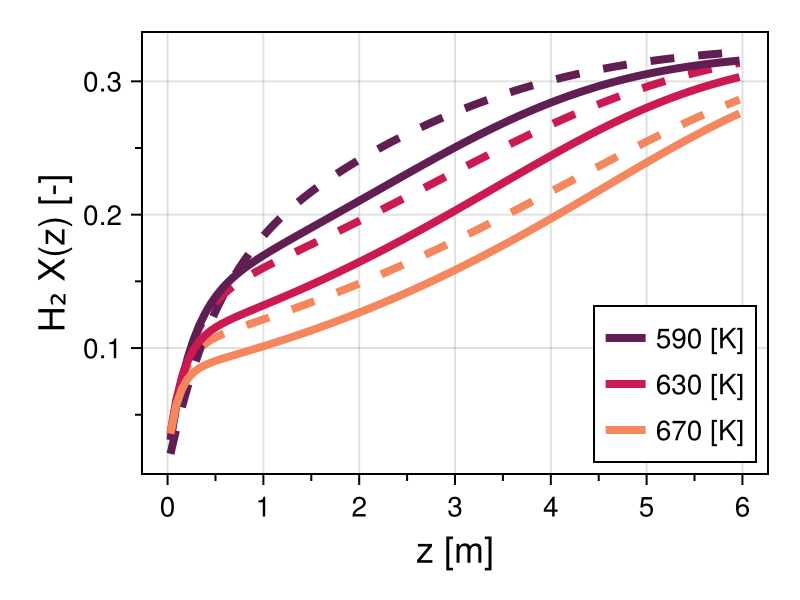}
        \caption{}
        \label{subfig:IDCRnominalconversionprofiles}
    \end{subfigure}
    \hfill
    \begin{subfigure}{0.325\textwidth}
        \includegraphics[width=\textwidth]{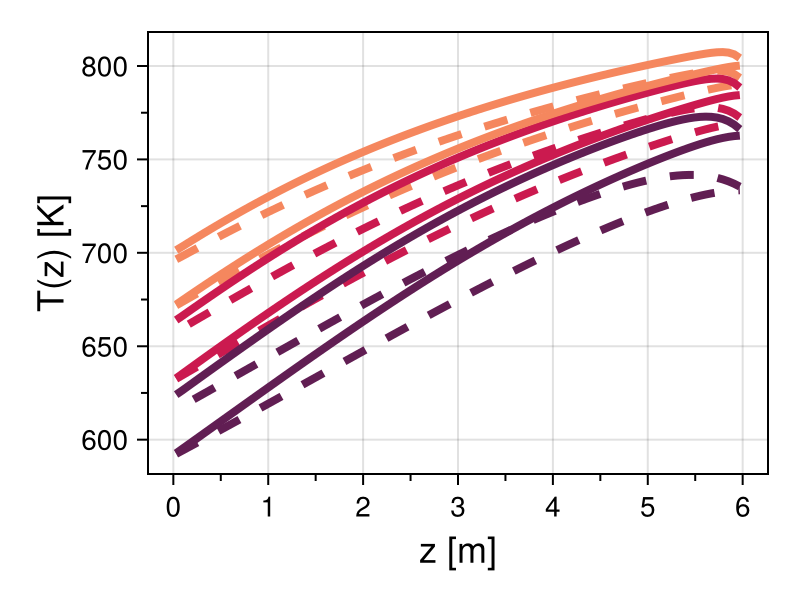}
        \caption{}
        \label{subfig:IDCRnominaltemperatureprofiles}
    \end{subfigure}
    \caption{IDCR steady states at nominal operating conditions. \textbf{(a):} Outlet H$_2$ conversion (yellow) and top temperature~(blue). \textbf{(b):} H$_2$ conversion profiles. \textbf{(c):} Temperature profiles. \textbf{(a)--(c):} Real-fluid properties (solid) and ideal-fluid properties (dashed).}
    \label{fig:IDCRnominal}
\end{figure*}

\noindent\textbf{IDCR}\\
\noindent 
Next, we perform simulations of the IDCR at nominal operating conditions. In Figure~\ref{fig:IDCRnominal}, we show the steady states for varying inlet temperatures $T_\text{in} \in [500,700]$. The steady-state curves in Figure~\ref{subfig:IDCRnominalcurves} depict the outlet H$_2$ conversion and the top temperature at the outlet of the IDCR's LCT. They display the characteristic S-shape, which indicates that the IDCR has multiple (and unstable) steady states in the interval of inlet temperatures between its extinction and ignition points \citep{rawlings2002chemical,khademi2017comparison,rosbo2024comparison}. Compared to the FBR, the IDCR achieves much higher H$_2$ conversion. The real model predicts an outlet conversion of $X_{\text{out}}^*=31.8\%$, which is slightly lower than the ideal models $x_{\text{out},\text{NH}_3}^*=32.4\%$. The higher conversion results from effective counter-current heat removal, which lowers the local temperature and thereby mitigates thermodynamic equilibrium limitations. This is clearly reflected by the ignited reactor profiles shown in Figure~\ref{subfig:IDCRnominalconversionprofiles}-\ref{subfig:IDCRnominaltemperatureprofiles}. At the lowest inlet temperature $T = 590$ [K], the ideal model operates near its optimum, where the balance between reaction kinetics and thermodynamic driving force is favorable. Hence, the conversion profile exhibits a consistently increasing slope, indicating a kinetic control throughout the reactor length. At higher inlet temperatures, both models reach their thermodynamic limits. This can be seen in the distinct kink in the conversion profiles within the first meter of the reactor.

\begin{figure}[tb]
    \centering
    \begin{subfigure}{\linewidth}
        \centering
        \includegraphics[width=0.75\linewidth]{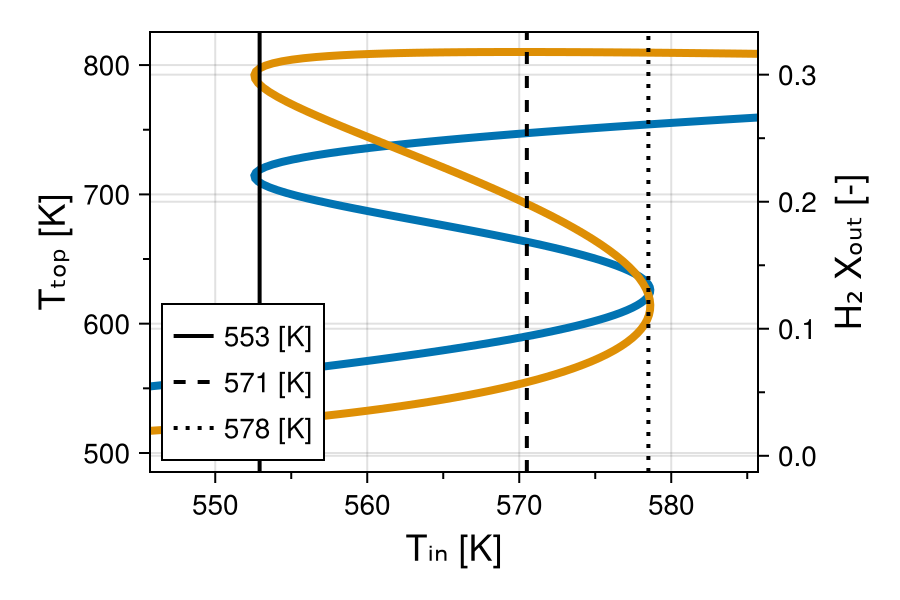}
        \caption{}
        \label{subfig:IDCRzoomreal}
    \end{subfigure}
    \begin{subfigure}{\linewidth}
        \centering
        \includegraphics[width=0.75\linewidth]{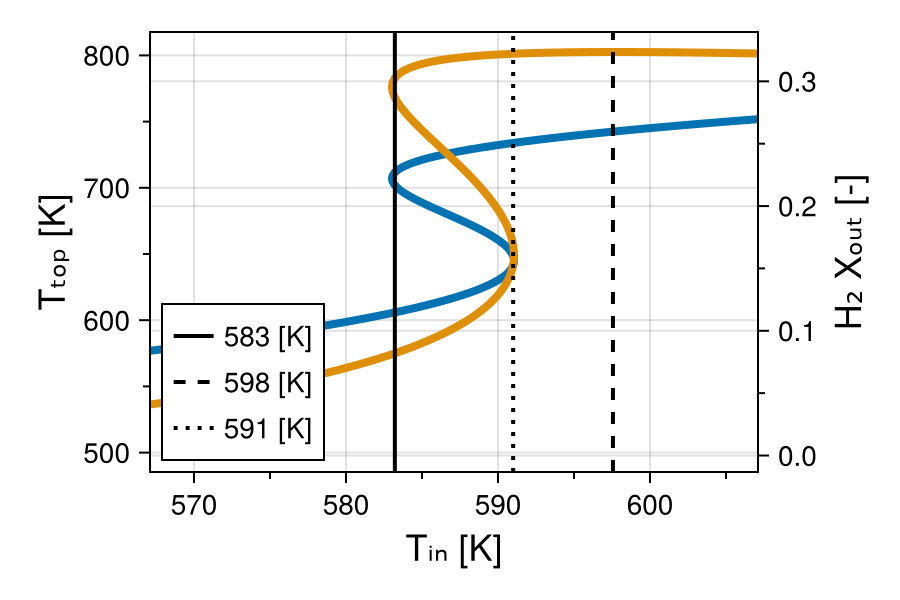}
        \caption{}
        \label{subfig:IDCRzoomideal}
    \end{subfigure}
    \caption{IDCR steady-state curves. Extinction, ignition, and optimal points. \textbf{(a):} Real fluid. \textbf{(b):} Ideal fluid.}
    \label{fig:IDCRzoom}
\end{figure}

\begin{figure}[tb]
    \centering
    \begin{subfigure}{\linewidth}
        \centering
        \includegraphics[width=0.75\linewidth]{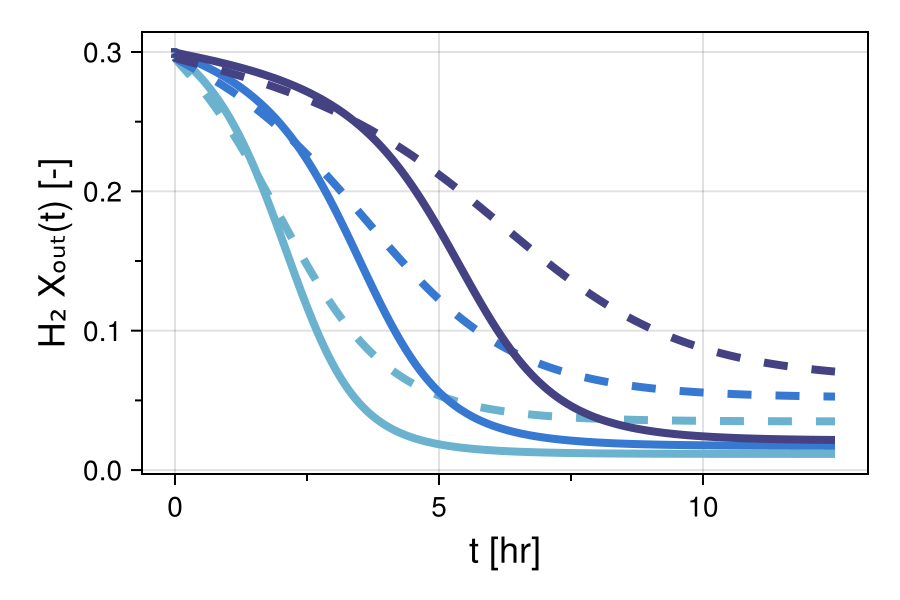}
        \caption{}
        \label{subfig:IDCRextinctionconversion}
    \end{subfigure}
    \begin{subfigure}{\linewidth}
        \centering
        \includegraphics[width=0.75\linewidth]{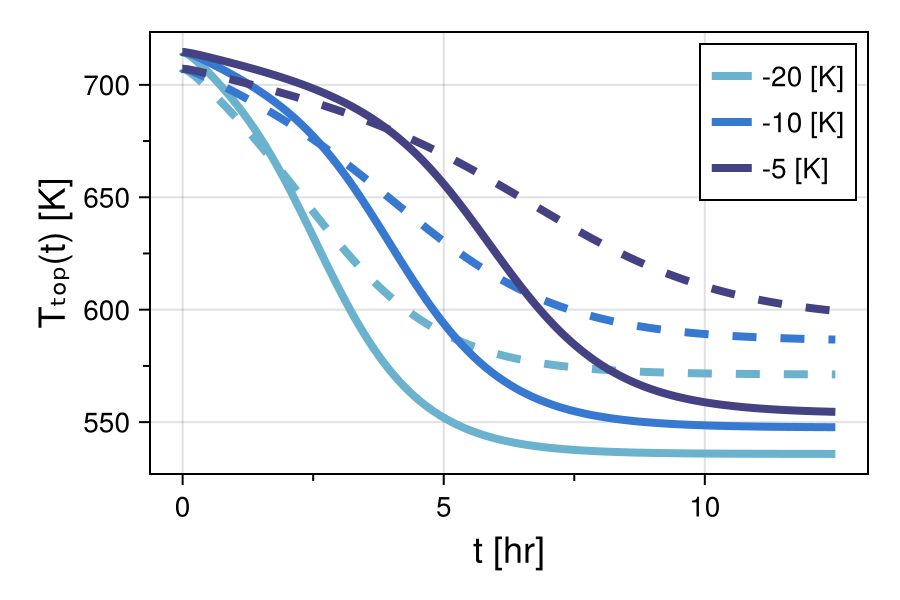}
        \caption{}
        \label{subfig:IDCRextinctiontemperature}
    \end{subfigure}
    \caption{IDCR simulations for varying steps in the inlet temperature away from the extinction point. \textbf{(a):} H$_2$ conversion. \textbf{(b):} Top temperature. \textbf{(a)--(b):} Real-fluid properties (solid) and ideal-fluid properties (dashed).}
    \label{fig:IDCRextinction}
\end{figure}

While both models exhibit qualitatively similar conversion profiles, there is a notable quantitative disparity between their S-shapes. In Figure~\ref{fig:IDCRzoom}, we depict the optimal inlet temperature together with the extinction and ignition points. The optimal inlet temperature of the real model is $T^*_\text{in}=571$ [K], which is 27 [K] lower than that of the ideal model at $T^*_\text{in}=598$ [K]. In addition, we see that the ideal model's extinction point lies 5 [K] higher than the real model's ignition point. Hence, a fully ignited steady-state for the real model may occur at operating conditions which cause an extinguished steady-state for the ideal model. Furthermore, we see that the interval of steady-state multiplicity for the ideal model spans 8 [K] which is contrasted by the three times larger 25 [K] interval for the real model. Hence, the choice of thermodynamic model may have a large implication on the chosen lower stability margin added to the optimal inlet temperature to prevent extinction, since it is safer to operate the IDCR too hot, rather than too cold. To illustrate the consequences of extinction, Figure~\ref{fig:IDCRextinction} shows dynamic simulations which step away from the ignited steady-state at the extinction point. Compared to the AFBR, we note that the dynamics happen on a time scale of hours, instead of a time scale of minutes. We note that the ideal model settles slightly later, and at higher conversion and top temperature. This is of course explained by its steady-state curve's more slender S-shape. 

\begin{figure}[tb]
    \centering
    \begin{subfigure}{\linewidth}
        \centering
        \includegraphics[width=0.75\linewidth]{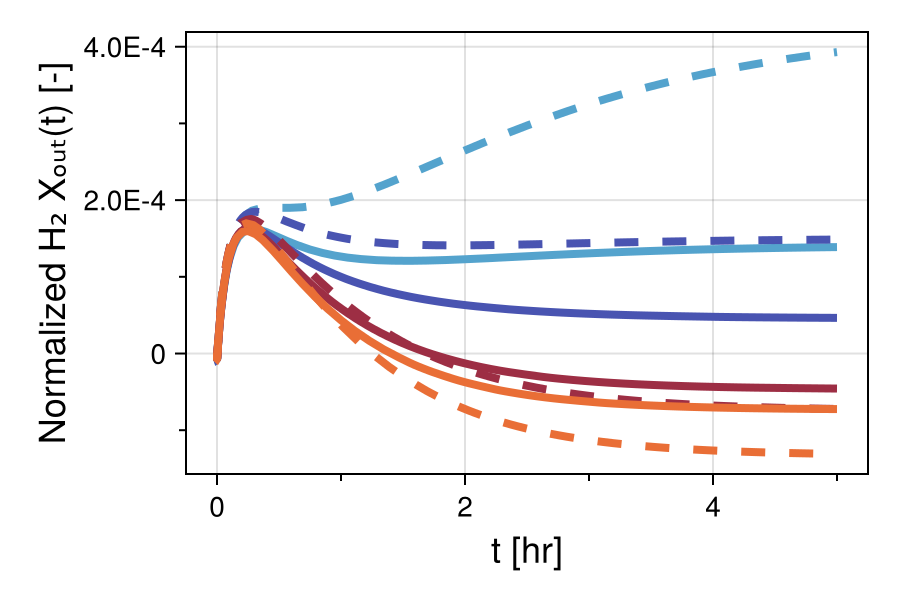}
        \caption{}
        \label{subfig:IDCRstepsconversion}
    \end{subfigure}
    \begin{subfigure}{\linewidth}
        \centering
        \includegraphics[width=0.75\linewidth]{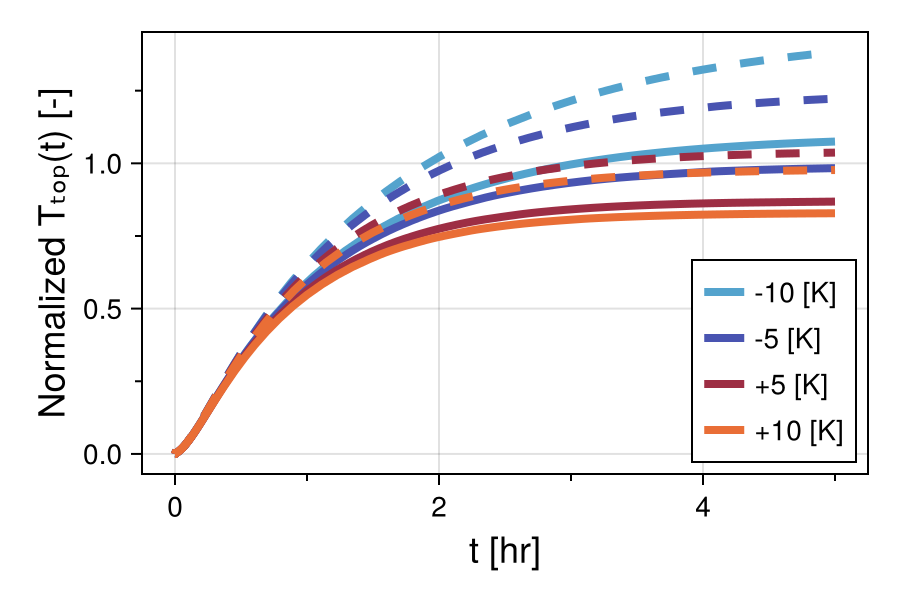}
        \caption{}
        \label{subfig:IDCRstepstemperature}
    \end{subfigure}
    \caption{IDCR step responses for varying steps in the inlet temperature at optimal operating conditions. \textbf{(a):} Normalized H$_2$ conversion. \textbf{(b):} Normalized top temperature. \textbf{(a)--(b):} Real-fluid properties (solid) and ideal-fluid properties (dashed).}
    \label{fig:IDCRsteps}
\end{figure}

For completeness, we also conduct step responses at the optimal operating conditions for the IDCR. They are shown in Figure~\ref{fig:IDCRsteps}. Similar consequences can be drawn. The ideal model settles more slowly, and is more sensitive to steps around its optimum.

\subsubsection{Literature model assumptions}
As discussed in Appendix~\ref{app:energybalance}, literature models apply ideal gas thermodynamics, neglect diffusive and conductive transport by setting $D = \kappa = 0$, and apply a pseudo-steady-state approximation by setting $\partial_t c = 0$. In our case study, these simplifications resulted in differences mostly within numerical integration tolerances, yielding simulation results nearly indistinguishable from those for the presented models. The most pronounced difference appeared in the short time-scale of the AFBR step-response simulations, as shown in Appendix~\ref{app:pseudosteps}.

%% file: tex/5conclusion.tex
\section{Conclusion}\label{sec:conclusion}
This work presented a flexible modeling framework for steady-state and dynamic simulation of fixed-bed reactor systems, demonstrated through an ammonia synthesis case study. The framework accommodates thermodynamic equations of state, consistent energy balances, and varying transport assumptions. Simulation results show that common literature assumptions, such as negligible dispersive transport and pseudo-steady-state mass balances, are valid for the presented case study. Moreover, while ideal thermodynamics suffice for simulating the AFBR, the slower dynamics of the IDCR reveal significant real-fluid effects, resulting in notably different steady-state behavior compared to ideal-fluid models. The results highlight the importance of selecting models that balance fidelity with simplicity and demonstrate the adaptability of the proposed framework to varying process configurations. The framework supports systematic model construction and can be extended beyond catalytic reactors to systems such as fixed-bed adsorption columns, membrane reactors, and separation processes. Future work will explore its use in reactor optimization and control for P2X applications.

%% file: tex/6appendix.tex
\appendix

\newcommand{\hbAppendixPrefix}{A.}
\renewcommand{\thefigure}{\hbAppendixPrefix\arabic{figure}}
\setcounter{figure}{0}
\renewcommand{\thetable}{\hbAppendixPrefix\arabic{table}} 
\setcounter{table}{0}
\renewcommand{\theequation}{\hbAppendixPrefix\arabic{equation}} 
\setcounter{equation}{0}

\section{Literature energy balance}\label{app:energybalance}
We present the energy balance in conservative form, derived directly in terms of the internal energy density. However, the energy balance is most often written in terms of temperature, and consequently, in non-conservative form. In this section, we will derive the latter form using common simplifying assumptions.

Using the definition of the fluid-phase internal energy density \eqref{eq:fluidenergydensity}, and applying the volume and internal energy constraints, \eqref{eq:volumeconstraint} and \eqref{eq:energyconstraint}, we find the relation
\begin{align}\label{eq:uhp}
    u=h-P,
\end{align}
where $h=H(T,P,c)$ [J$/$m$^3$] denotes the enthalpy density. As a consequence of degree-one homogeneity of the thermodynamic functions, we may apply Euler's first theorem for homogeneous functions \citep{michelsen2008thermodynamic}, and find that the enthalpy density can be written in terms of the partial molar enthalpies $\bar{H}=[\bar{H}_\alpha]_{\alpha \in \mathcal{C}}$~[J$/$mol], i.e.,
\begin{align}
    h=H(T,P,c)=\sum_{\alpha \in \mathcal{C}}c_\alpha \bar{H}_\alpha.
\end{align}
Assuming an ideal fluid, the partial molar enthalpies correspond to those given by the pure component properties, independent of pressure and composition. Consequently, we find the total differential
\begin{align}
    \partial h=\sum_{\alpha \in \mathcal{C}}\left(\partial c_\alpha \bar{H}_\alpha + c_\alpha\partial_T \bar H_\alpha \partial T\right).
\end{align}
By definition of the component-wise specific heat capacity $c_{p,\alpha}$ [J$/$(kg$\cdot$K)] and molecular weight $M_\alpha$ [kg$/$mol], we get
\begin{align}
    \partial_T\bar{H}_\alpha=M_\alpha c_{p,\alpha}.
\end{align}
Invoking the definition of fluid density \eqref{eq:fluiddensity}, we then find that
\begin{align}\label{eq:dh}
    \partial h=\sum_{\alpha \in \mathcal{C}}\partial c_\alpha \bar{H}_\alpha + \rho c_{p} \partial T,
\end{align}
where $c_{p}=\sum_{\alpha \in \mathcal{C}}\omega_\alpha c_{p,\alpha}$ [J$/$(kg$\cdot$K)] denotes the mass-averaged specific heat capacity, specified by the mass fractions $\omega_\alpha=M_\alpha c_\alpha /\rho$ [--]. If we neglect dispersive effects, the molar flux reduces to the advective flux $N=v\odot c$ and the energy flux reduces to the convective flux $E=H(T,P,N)$. Moreover, we similarly find the total differential
\begin{align}\label{eq:de}
    \partial E=\sum_{\alpha \in \mathcal{C}}\partial N_\alpha \bar{H}_\alpha + \rho c_{p}v \partial T,
\end{align}
which is written in non-conservative form.

Putting it all together, we first insert \eqref{eq:uhp} into the energy balance \eqref{eq:energybalances}, and find the enthalpy balance
\begin{align}\label{eq:enthalpybalance}
    \partial_th=-\partial_zE+Q,
\end{align}
assuming constant pressure in time ($\partial_t P = 0$) . Rearranging the component-wise mass balance of \eqref{eq:massbalances} yields
\begin{align}
    \partial_tc_\alpha+\partial_zN_\alpha=R_\alpha,
\end{align}
which, combined with both \eqref{eq:dh}-\eqref{eq:de}, we use to expand the enthalpy balance \eqref{eq:enthalpybalance}
\begin{align}\label{eq:firststep}
    \rho c_{p} \left(\partial_t T+v\partial_zT\right)=-\sum_{\alpha \in \mathcal{C}} R_\alpha \bar{H}_\alpha + Q.
\end{align}
Inserting the definition of the component-wise reaction rates~\eqref{eq:compreaction}, and exchanging order of summation, we find that
\begin{align}\label{eq:reactionheat}
    \sum_{\alpha \in \mathcal{C}}R_\alpha \bar{H}_\alpha = \sum_{k\in\mathcal{R}}r_k\Delta H_k,
\end{align}
where we have defined the heats of reaction $\Delta H=[\Delta H_k]_{k\in\mathcal{R}}$ [J/mol] by the reaction stoichiometry
\begin{align}\label{eq:heatofreaction}
    \Delta H_k=\sum_{\alpha\in\mathcal{C}}\nu_{k,\alpha}\bar{H}_\alpha.
\end{align}
In summary, inserting \eqref{eq:reactionheat} into \eqref{eq:firststep} gives the balance
\begin{align}
    \rho c_p\partial_tT=-\rho c_p v \partial_zT-\Delta H^Tr + Q.
\end{align}
Extending the discussion to heterogeneous volumes, we find the pseudo-homogeneous total energy balance
\begin{align}\label{eq:litenergybal}
\begin{split}
    \big(\varepsilon&\rho_\text{fluid} c_{p,\text{fluid}}+(1-\varepsilon)\rho_\text{solid}c_{p,\text{solid}}\big) \partial_tT\\
    &=-\varepsilon\left( \rho_\text{fluid} c_{p,\text{fluid}} v\partial_z T -\Delta H^T r\right) +Q,
\end{split}
\end{align}
where $c_{p,\text{fluid}},c_{p,\text{solid}}$ [J/(kg$\cdot$K)] and $\rho_\text{fluid},\rho_\text{solid}$ [kg/m$^3$] denote the specific heat capacities and densities of the fluid and solid phases. Equation \eqref{eq:litenergybal} corresponds to the formulation typically found in literature, which is often accompanied by the pseudo-steady-state mass balance equations
\begin{align}\label{eq:pseudo-steady-state}
    \partial_zN=R.
\end{align}

\section{Step responses with literature model assumptions}\label{app:pseudosteps}
\begin{figure}[tb]
    \centering
    \includegraphics[width=0.75\linewidth]{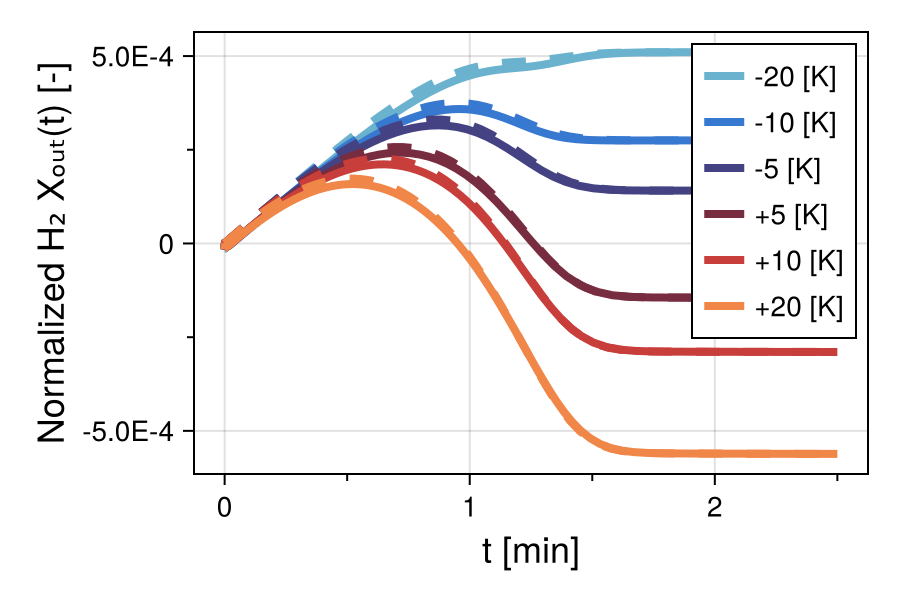}
    \caption{AFBR normalized H$_2$ conversion step responses for varying steps in the inlet temperature at optimal operating conditions. Presented model assumptions (solid) and literature model assumptions (dashed).}
    \label{fig:AFBRliteratureassumptions}
\end{figure}
Figure~\ref{fig:AFBRliteratureassumptions} shows the H$_2$ conversion step responses for the presented model and for the same model using the literature assumptions described in Appendix~\ref{app:energybalance}. That is, we assume ideal-gas thermodynamics, negligible dispersive terms, and pseudo-steady-state mass balances. The differences are negligible, with the pseudo-steady-state assumption only causing a barely noticeable deviation during the transient.